\documentclass[12pt]{amsart}

\usepackage{amssymb,amsfonts}
\usepackage{amsthm}
\usepackage{graphicx}
\usepackage[all,arc]{xy}
\usepackage[colorlinks=true, allcolors=blue]{hyperref}
\usepackage{enumerate}
\usepackage{bbm}
\usepackage{dsfont}
\usepackage{mathrsfs}
\usepackage{tikz-cd}
\usetikzlibrary{shapes}
\usepackage{import}
\usepackage{float}
\usepackage{pstricks}
\usepackage{pst-plot}
\usepackage{url}
\restylefloat{table}

\usepackage[left=1in,top=1in,right=1in,bottom=1in]{geometry}

\usepackage{mathtools}

\usepackage[normalem]{ulem}   

\newtheorem{myth}{Theorem}[section]
\newtheorem{mycor}[myth]{Corollary}
\newtheorem{myprop}[myth]{Proposition}
\newtheorem{mylemma}[myth]{Lemma}

\newtheorem{theorem}{Theorem} 

\newtheorem*{theorem*}{Theorem A} 
\newtheorem*{theoremb*}{Theorem B} 
\newtheorem*{theoremc*}{Theorem C}

\theoremstyle{definition}
\newtheorem{mydef}[myth]{Definition}

\newtheorem{myex}[myth]{Example}

\newtheorem{myrem}[myth]{Remark}

\theoremstyle{remark}

\newtheorem{sch}[myth]{Scholium}

\newcommand{\Z}{\mathbb{Z}}

\newcommand{\R}{\mathbf{R}}

\renewcommand{\emptyset}{\varnothing}

\newcommand{\G}{\mathbf{G}}

\usepackage[colorinlistoftodos]{todonotes}

\def\G{\mathbf{G}}

\def\R{\mathbb{R}}

\def\Z{\mathbb{Z}}

\newcommand*\circled[1]{\tikz[baseline=(char.base)]{\node[shape=circle,draw,inner sep=2pt, scale=0.7] (char) {#1};}}
\newcommand{{\circone}}{\circled{$1$}}
\newcommand{{\circtwo}}{\circled{$2$}}
\newcommand{{\circthree}}{\circled{$3$}}
\newcommand{{\circn}}{\circled{$n$}}

\newtheorem*{thm*}{Theorem}

\makeatletter

\title{Filtrations of MODULI SPACES OF TROPICAL  WEIGHTED  STABLE CURVES}

\author[S. Serpente]{Stefano Serpente}\address{Dipartimento di Matematica e Fisica, Università Roma Tre, Roma, I-00146}\email{\url{stefano.serpente@uniroma3.it}}

\begin{document}
\begin{abstract}
We consider tropical versions of Hassett's moduli spaces of weighted stable curves $M_{g,\mathcal{A}}^{trop}$, $\overline{M}_{g,\mathcal{A}}^{trop}$ and $\Delta_{g,\mathcal{A}}$ associated to a weight datum $\mathcal{A}=(a_1,...,a_n)\in(\mathbb{Q}\cap(0,1])^n$, their associated graph complexes $G^{(g,\mathcal{A})}$ and study the topology of these spaces as $\mathcal{A}$ changes.  We show that for fixed $g$ and $n$, there are particular filtrations of these topological spaces and their graph complexes which may be used to compute the reduced rational homology of $\Delta_{g,\mathcal{A}}$ and the top weight cohomology of the moduli space $\mathcal{M}_{g,\mathcal{A}}$ of smooth $(g,\mathcal{A})$-stable algebraic curves.
\end{abstract}
\maketitle

\tableofcontents

\section{Introduction}
 Let $g\geq 0$ and $n\geq 1$ be two integers, and $\mathcal{A}=(a_1,...,a_n)\in\mathcal{D}_{g,n}:=(\mathbb{Q}\cap(0,1])^n$ a weight datum such that $2g-2+\sum_{i=1}^n a_i>1.$  In order to give different compactifications to the moduli stack of smooth curves $\mathcal{M}_{g,n}$, Hassett in \cite{HASSETT2003316} defined generalized stability conditions on algebraic curves introducing weight data attached to marked points. The definition of weighted stable curve gave rise to the moduli stacks $\overline{\mathcal{M}}_{g,\mathcal{A}}$, which generalize the standard moduli stack of stable curves $\overline{\mathcal{M}}_{g,n}$. In \cite{ulirsch14tropical}, Ulirsch exploited the dualism between the moduli theory of algebraic curves and the one of tropical curves defining generalized stability conditions on graphs, and out of these, constructed the moduli spaces of tropical weighted stable curves $M_{g,\mathcal{A}}^{trop}$, generalizing the construction of the moduli spaces of tropical curves with marked points $M_{g,n}^{trop}$ from \cite{article} (see Section \ref{section2}).

Following \cite{HASSETT2003316}, the space of admissible weight data $\mathcal{D}_{g,n}$ may be decomposed into connected components separated by a set of hyperplanes, which are called walls. The set of these hyperplanes is called chamber decomposition, and the connected components of the decomposition are called chambers. We call $\mathbf{K}$ the set of chambers. We treat the wall crossing theory for these stability conditions in Section \ref{section3}. Each wall can be defined by a linear inequality and subdivides the space of weights into two chambers, according to the two corresponding inequality, see for example \ref{Dg2} and \ref{Dg3}.

We consider the fine chamber decomposition, defined in \cite{HASSETT2003316}:$$\mathit{W}_f=\left\lbrace \sum_{j\in S}a_j=1:\lbrace S\subset \lbrace 1,...,n \rbrace\rbrace, 2\leq|S|\leq n-2\delta_{g,0} \right\rbrace,$$
where $\delta_{i,j}$ is the Kronecker Delta. As showed by Ulirsch in \cite{ulirsch14tropical}, different weight data in the same chamber give rise to the same moduli space, and the fine chamber decomposition is the coarsest one with this property. 

Given two  weight data $\mathcal{A}$ and $\mathcal{B}$, we write $\mathcal{A}\leq \mathcal{B}$ if $a_i\leq b_i$ component-wise. It is easy to show that if $\mathcal{A}\leq \mathcal{B}$ there is an inclusion as a closed subset

\begin{equation}\label{one}
   M_{g,\mathcal{A}}^{trop} \subseteq M_{g,\mathcal{B}}^{trop}.
\end{equation} 
 
 For a given weight datum $\mathcal{A}$, we denote by $Ch_{\mathcal{A}}$ its chamber. We define a partial order relation in the set $\mathbf{K}$ which extends the previous one on weight data as follows. Let $Ch_1,Ch_2\in\mathbf{K}$, we say that $Ch_1\leq Ch_2$ if they are equal or there are $S_1,...,S_t\subset \lbrace 1,...,n\rbrace$ such that for every $\mathcal{A}\in Ch_1$, $\sum_{i\in S_j}a_i<1$ and  for every $\mathcal{B}\in Ch_2$, $\sum_{i\in S_j}a_i>1$, for every $j$ from 1 to $t$, while for any $S'\neq S_j$ for every $j$ from 1 to $t$, the two chambers belong to the same half-space of $\mathcal{D}_{g,n}$ induced by the wall $\sum_{i\in S'}a_i=1$. It is easy to see that if we pick $\mathcal{A}\in Ch_1$ and $\mathcal{B}\in Ch_2$ with $Ch_1\leq Ch_2$, the inclusion (\ref{one}) still holds. 

 Given a weight datum $\mathcal{A}=(a_1,...,a_n)$ and a permutation $\sigma \in S_n$, let $\sigma (\mathcal{A})$ be the weight datum obtained by permuting the weights of $\mathcal{A}$ according to $\sigma$, i.e. $\sigma(\mathcal{A})=(a_{\sigma (1)},...,a_{\sigma(n)})$. Then $M_{g,\mathcal{A}}^{trop}$  is homeomorphic to $M_{g,\sigma (\mathcal{A})}^{trop}$: the homeomorphism, called relabeling homeomorphism, consists of sending a tropical curve into the curve with the same underlying graph and the same length function, but with the legs marked according to the permutation (see e.g. \ref{exSn}).

 We can consider the action of $S_n$ induced on $\mathbf{K}$ given by $\sigma(Ch_{\mathcal{A}})=Ch_{\sigma(\mathcal{A})}$ as well. We call the orbits of this action chambers up to symmetry, and we denote the set of chambers up to symmetry by $[\mathbf{K}]$. Denote by $[Ch]$ the orbit of $Ch\in\mathbf{K}$, which we call the chamber up to symmetry of $Ch$. We say that $[Ch_1]\leq[Ch_2]$ if there are two chambers $Ch_{\mathcal{A}_1}\in [Ch_1]$ and $Ch_{\mathcal{A}_2}\in[Ch_2]$ such that $Ch_{\mathcal{A}_1}\leq Ch_{\mathcal{A}_2}$, for some $\mathcal{A}_1$ and $\mathcal{A}_2$. In that case, each time we pick $$\mathcal{A}\in Ch_{\mathcal{A}_1}\in [Ch_{\mathcal{A}_1}]\leq[Ch_{\mathcal{A}_2}]\ni Ch_{\mathcal{A}_2}\ni\mathcal{B},$$ there is a permutation $\sigma\in S_n$ giving a topological embedding $$M_{g,\mathcal{A}}^{trop}\hookrightarrow M_{g,\mathcal{B}}^{trop}$$ 
obtained combining the relabeling homeomorphism $M_{g,\mathcal{A}}^{trop}\cong M_{g,\sigma(\mathcal{A})}^{trop}$ with the inclusion (\ref{one}) $M_{g,\sigma(\mathcal{A})}^{trop}\subset M_{g,\mathcal{B}}^{trop}$. In particular, there is only a moduli space of tropical curves up to homeomorphism for each chamber up to symmetry.

 All these properties work if we replace $M_{g,\mathcal{A}}^{trop}$ with $\overline{M}_{g,\mathcal{A}}^{trop}$ and $\Delta_{g,\mathcal{A}}$, respectively the moduli space of extended $\mathcal{A}$-weighted stable curves of genus $g$ and the moduli space of $\mathcal{A}$-weighted stable curves of genus $g$ and volume 1 (see Section \ref{moduli} for their definitions).

 We illustrate the situation with an example. Suppose we have $g\geq 1$ and $n=3$, and we have the weight data $\mathcal{A}_1=(\frac{12}{27},\frac{14}{27},1-\varepsilon)$ and $\mathcal{A}_2=(\frac{14}{27}-\varepsilon, \frac{12}{27}, \frac{14}{27})$, for some $0< \varepsilon <\frac{1}{27}$. Clearly they are not comparable with respect to the partial order on weight data, and we can verify that they belong to different chambers, i.e. $Ch_{\mathcal{A}_1}\neq Ch_{\mathcal{A}_2}$. It is also possible to verify that $[Ch_{\mathcal{A}_1}]\neq [Ch_{\mathcal{A}_2}]$. This implies that their moduli spaces $M_{g,\mathcal{A}_1}^{trop}$ and $M_{g,\mathcal{A}_2}^{trop}$ are different, and none of them is the subspace of the other. But if we reorder the weights of the second
weight by the permutation $\sigma=(1\quad 3\quad  2 )\in S_3 $ the datum we obtain is $\sigma(\mathcal{A}_2)=(\frac{12}{27},\frac{14}{27}, \frac{14}{27}-\varepsilon)$, and by the relabeling homemorphism we know that
$$M_{g,\mathcal{A}_2}^{trop}\approx  M_{g,\sigma(\mathcal{A}_2)}^{trop}.$$

Moreover $\sigma(\mathcal{A}_2)$ and $\mathcal{A}_1$ are comparable, in particular $\sigma(\mathcal{A}_2)\leq \mathcal{A}_1$, so we have the inclusion as a subspace $$M_{g,\sigma(\mathcal{A}_2)}^{trop} \subset M_{g,\mathcal{A}_1}^{trop}.$$ 
Composing the relabeling homeomorphism with this inclusion gives an embedding as a subspace $$M_{g,\mathcal{A}_2}^{trop} \hookrightarrow M_{g,\mathcal{A}_1}^{trop}.$$ 
 In the case $g\geq 1$ and $n=3$ there are five chambers up to symmetry, so choosing representative weight data up to relabeling homeomorphisms we get the filtration $$ M_{g,(\frac{1}{3},\frac{1}{3},\frac{1}{3}-\varepsilon)}^{trop} \subset  M_{g,(\frac{4}{9}-\varepsilon,\frac{4}{9}-\varepsilon,\frac{4}{9}-\varepsilon)}^{trop} \subset M_{g,(\frac{12}{27},\frac{14}{27}, \frac{14}{27}-\varepsilon)}^{trop} \subset M_{g,(\frac{12}{27},\frac{14}{27},1-\varepsilon)}^{trop} \subset M_{g,3}^{trop}$$ where we take $\varepsilon$'s in order to take weight data in the interior of the chamber decomposition. 
 (see Example \ref{final} for further details). 

 Following the previous discussion, in Section \ref{section4}, we prove the following Theorem.
Let $g\geq 0$, $n\geq 1$ and $\mathcal{A}\in\mathcal{D}_{g,n}$ $C\subset M_{g,\mathcal{A}}^{trop}$ be a closed subset. We say it is a sub-moduli space if it is homeomorphic to $M_{g,\mathcal{B}}^{trop}$ for some $\mathcal{B}$, and points of $C$ are in bijection with $(g,\mathcal{B})$-stable tropical curves.
\begin{theorem}\label{Main1}
Let $g\geq 0$, $n\geq 1$ be two integers. Fix a weight datum $\mathcal{A}\in \mathcal{D}_{g,n}$.
There are filtrations of $M^{trop}_{g,\mathcal{A}}$ given by embeddings as sub-moduli spaces induced by the partial order on the set of chambers up to symmetry. Namely, given a sequence  $$[Ch_{\mathcal{A}_1}]\leq [Ch_{\mathcal{A}_2}]\leq ...\leq [Ch_{\mathcal{A}_p}] \leq ... \leq [Ch_{\mathcal{A}_{N-1}}]\leq [Ch_{\mathcal{A}}], $$ the filtration is
    $$M^{trop}_{g,\mathcal{\mathcal{A}}_1}\hookrightarrow M^{trop}_{g,\mathcal{A}_2}\hookrightarrow  ...\hookrightarrow M^{trop}_{g,\mathcal{A}_p} \hookrightarrow ... \hookrightarrow M^{trop}_{g,\mathcal{A}_{N-1}}\hookrightarrow M^{trop}_{g,\mathcal{A}}. $$
 The same result holds if we replace $M_{g,\mathcal{A}}^{trop}$ with the moduli space of extended weighted tropical $(g,\mathcal{A})$-stable curves $\overline{M}_{g,\mathcal{A}}^{trop}$  or the moduli space of $(g,\mathcal{A})$-stable tropical curves of volume 1 $\Delta_{g,\mathcal{A}}$.
\end{theorem}

In Section \ref{section4}, we define the graph complexes $G^{(g,\mathcal{A})}$. These objects generalize the graph complexes $G^{(g,n)}$ defined in section 2.4 of \cite{CGP2}. We deduce a Filtration Theorem analogous to Theorem \ref{Main1} for these complexes (see Section \ref{section4} for further details). This gives them the structure of filtered chain complexes

 In order to study further these graph complexes and their homology, in Section \ref{section6} we extend to $\Delta_{g,\mathcal{A}}$ the theory developed in \cite{CGP1} and \cite{CGP2} for $\Delta_{g,n}$, treating it as a symmetric $\Delta$-complex for every $\mathcal{A}$, and then we generalize Theorem 1.4 of \cite{CGP2} showing that there is a natural surjection of chain complexes $C_*(\Delta_{g,\mathcal{A}})\rightarrow G^{(g,\mathcal{A})}$ decreasing degrees by $2g-1$ inducing isomorphisms on homology $$\widetilde{H}_{k+2g-1}(\Delta_{g,\mathcal{A}};\mathbb{Q})\rightarrow H_k(G^{(g,\mathcal{A})})$$ for all $k$'s.  Analogously to what happen for $\Delta_{g,n}$
there is a natural isomorphism $$Gr^W_{6g-6+2n}H^{4g-6+2n-k}(\mathcal{M}_{g,\mathcal{A}};\mathbb{Q})\rightarrow {H}_{k}(G^{(g,\mathcal{A})}),$$ between the rational top weight cohomology of $\mathcal{M}_{g,\mathcal{A}}$ and the rational homology of the complex $G^{(g,\mathcal{A})}$.

There is a spectral sequence associated to a filtered chain complex which can be used to compute the homology of the complex. In particular, the structure of bounded filtered chain complex given to each $G^{(g,\mathcal{A})}$, combined with the shifting degree isomorphism of the top weight cohomology of $\mathcal{M}_{g,\mathcal{A}}$ with the homology of the complex gives us the following Theorem, proved at the end of the section \ref{section6}.
\begin{theorem}\label{LastTheorem}
Fix $g\geq 1$, $n\geq 2$. Assume we have a sequence of chambers up to symmetry $[Ch_{\mathcal{A}_1}]\leq...\leq [Ch_{\mathcal{A}_p}]\leq...\leq [Ch_{\mathcal{A}_N}]$, and let $G^{(g,\mathcal{A}_1)}\hookrightarrow...\hookrightarrow G^{(g,\mathcal{A}_p)}\hookrightarrow...\hookrightarrow G^{(g,\mathcal{A}_N)}$ be the induced filtration on graph complexes. Then 
$$Gr^W_{6g-6+2n}H^{4g-6+2n-k}(\mathcal{M}_{g,\mathcal{A}_N};\mathbb{Q})\cong \bigoplus_{p=1}^NE^{\infty}_{p,k-p},$$
where the terms $E^{\infty}_{p,k-p}$ are the ones to which the spectral sequence induced by the above filtration converges.
\end{theorem}

In Example \ref{FtComp}, we use this Theorem to compute the Top Weight Cohomology of $\mathcal{M}_{1,3}$, confirming results obtained in \cite{CGP2}.

\subsection{Motivation and related works}
This work builds up on previous work and ideas from many authors in the area. Firstly, the moduli space of weighted tropical curves was constructed by Ulirsch in \cite{ulirsch14tropical}, Propositions 6.1 and 6.2 where the author shows the constancy of the tropical moduli spaces inside the chambers. In the same work, the moduli spaces of weighted tropical curves are identified with the skeletons of moduli spaces of curves, generalizing \cite{ACP} in the case of moduli spaces of curves with standard stability. The inclusion property coming from the componentwise relation on weight data was known in other works related to tropical moduli spaces such as \cite{CMPRS}. The relabeling symmetry in the algebraic set was a well known fact since Hassett published his  work \cite{HASSETT2003316}, and again the permutation action was studied in \cite{CMPRS} on $\Delta_{0,w}$ for certain particular cases, in \cite{kannan2021symmetries} to study $Aut(\Delta_{g,n})$ and in \cite{CFGP} for computing the $S_n$-equivariant cohomology of $\mathcal{M}_{g,n}$. In \cite{yun2020sn}, there are computations for the $S_n$-equivariant rational homology of the tropical moduli spaces $\Delta_{2,n}$ for $n\leq 8$. In \cite{hahn2021intersecting} authors give an explicit formula for the intersection products of weighted tropical $\psi$-classes on $M^{trop}_{0,\mathcal{A}}$, in arbitrary
dimensions. Sections 5 and 6 develop on tools and language introduced in \cite{CGP1} and \cite{CGP2}, which we adapt to study $\Delta_{g,\mathcal{A}}$.

Furthermore, there are many works on the topology of moduli spaces of tropical curves and their identifications with other objects from which we benefit:
\begin{itemize}

    \item In \cite{allcock}, Allcock, Corey and Payne showed that $\Delta_{g,n}$ is simply connected for $(g,n)\neq(0,4),(0,5)$. This result is generalized for $g\geq 1$ and for all $\mathcal{A}$ in \cite{KLSY} by Kannan, Li, Yun and the author.
    \item In \cite{CGP1} and \cite{CGP2} Chan, Galatius and Payne give results on the top weight cohomology of $\Delta_{g,n}$. In \cite{KLSY} similar results are given for $\Delta_{g,\mathcal{A}}$.
    \item When $\mathcal{A} = 1^{(n)}$, Vogtmann showed that $\Delta_{0, n}$ is homotopic to a wedge of $(n - 2)!$ spheres of dimension $n-4$, see \cite{ vogtmann_1990}.
    \item When $\mathcal{A}$ is heavy/light, i.e. it has $m$ components equal to $\varepsilon< 1/m$ and $n-m$ components equal to 1, Cavalieri, Hampe, Markwig, and Ranganathan in \cite{CHMR2014moduli}  derived from their results that $\Delta_{0, w}$ is homotopic to a wedge of $(n-2)!(n-1)^{m}$ spheres of dimension $n + m - 4$. 
    \item When $w$ has at least two weight-$1$ entries, Cerbu, Marcus, Peilen, Ranganathan and  Salmon in \cite{CMPRS} showed that $\Delta_{0, w}$ is homotopic to a wedge of spheres of possibly varying dimensions, and gave infinite families of $w$ where $\Delta_{0, w}$ is disconnected, and examples where $\pi_{1}(\Delta_{0, w}) = \Z/2\Z$.  
    \item In  \cite{freedman2021automorphisms} they study the automorphism group of $\Delta_{g,\mathcal{A}}$, which correspond to the subgroup of the group of relabeling morphisms preserving $\Delta_{g,\mathcal{A}}$.
    \end{itemize}
\subsection{Acknowledgements}  I am very grateful to my Ph.D. supervisor Margarida Melo for helping me through the writing of this work both from the mathematical point of view and in the form. I want also like to thank Melody Chan, Siddarth Kannan, Shiyue Li, Martin Ulirsch and Claudia Yun for the many discussions about the topic and this work. Thanks to the anonymous referee for many important suggestions and comments.

 \section{Background}\label{section2}
\subsection{Graphs} We start by introducing all the notation and the conventions we will use during the rest of the paper.
\begin{mydef} A decorated graph with $n$ legs $G$ is the data of:
\begin{itemize}
\item[1)] A finite non-empty set $V(G)$ called the set of vertices;
\item[2)] A finite set of half-edges $H(G)$;
\item[3)] An involution $\iota:H(G)\rightarrow H(G)$ with $n$ fixed elements, called legs, whose set is denoted by $L(G)$;
\item[4)] An endpoint map $\epsilon:H(G)\rightarrow V(G)$;
\item[5)] A vertex weight function on the vertices $w:V(G)\rightarrow \mathbb{Z}_{\geq 0}$.
\end{itemize}
\end{mydef}
A non-ordered pair $e=\lbrace h, h'\rbrace$ of distinct elements in $H(G)$ interchanged by the involution is an edge of the graph, and the set of edges is denoted by $E(G)$. If $\epsilon(h)=v$ we say that $h$ is adjacent to $v$, and that $v$ is the endpoint of $h$. The same definition works for edges. The valence of a vertex $v$ is $\text{val}(v):=|\epsilon^{-1}(v)|$, i.e., the number of half-edges adjacent to $v$. An edge whose endpoints coincide is called a loop. Two different edges with the same endpoints are said to be parallel. Two legs are called disjoint if their endpoints are distinct.

In the literature, decorated graphs are known also as weighted graphs, and the vertex weight function is known also as weight function. Here we change the terminology in order to avoid confusion with another notion of weight we will introduce later.
\begin{mydef}
The genus of a decorated graph is
\begin{center}
$g(G)=b_1(G)+\sum\limits_{v\in V(G)} w(v)$\\
\end{center}
whit $b_1(G):=|E(G)|-|V(G)|+c,$
where $c$ is the number of connected components of the graph.
\end{mydef}
From now on we will consider only connected graphs, so $c=1$ every time, and we will omit the adjective connected. 
\begin{mydef}
A morphism between graphs $G$ and $G'$ is a map
$$\alpha:V(G)\cup H(G)\rightarrow V(G')\cup H(G')$$
such that $\alpha (L(G))\subset L(G')$ and the following diagrams commute:
\begin{center}
\begin{tikzcd}[column sep=small]
V(G)\cup H(G) \arrow[r,"\alpha "] \arrow[d,"(id_{V\prime}  \epsilon)"]
& V(G')\cup H(G') \arrow[d,"(id_{V'\prime} \epsilon ')"]
&
&V(G)\cup H(G) \arrow[r,"\alpha "] \arrow[d,"(id_{V\prime} \iota)"]
& V(G')\cup H(G') \arrow[d,"(id_{V'\prime } \iota ')"]\\
V(G)\cup H(G) \arrow[r,"\alpha "]
& V(G')\cup H(G')
&
&V(G)\cup H(G) \arrow[r,"\alpha "]
& V(G')\cup H(G').
\end{tikzcd}
\end{center}
\end{mydef}
In particular, we have that $\alpha(V(G))\subseteq V(G')$.

A morphism $\alpha:G\rightarrow G'$ is said to be an isomorphism if it induces by restriction three bijections: $\alpha_V:V(G)\rightarrow V(G'),$ 
$\alpha_E:E(G)\rightarrow E(G')$ and
$\alpha_L:L(G)\rightarrow L(G')$. 
An automorphism of a graph $G$ is an isomorphism of $G$ with itself: in this case  endpoints of the legs are preserved.

 If the image of an edge $e$ is $v'\in V(G')$,  also its endpoints are mapped into $v'$ and we say that $e$ is contracted by $\alpha$ and that $\alpha$ is a contraction. Let $T\subseteq E(G)$, the graph $G/T$ denotes the graph obtained by contracting the edges $e\in T$. A weighted contraction is $(G/T, w/T)$ where $G/T$ is a contraction and $w/T$ is the vertex weight function defined by setting, for every $v\in V(G/T)$ 
\begin{equation}\label{contreq}
w/T(v)=b_1(\alpha^{-1}(v))+\sum\limits_{u\in \alpha^{-1}(v)} w(u).
\end{equation}  
 It is clear from equation (\ref{contreq}) that the genus of a decorated graph remains constant after contraction.

 Let $G$ be a decorated graph, and let $L(G)$ be the set of its legs. A marking is the assignment of a number from 1 to $n$ to each leg, and a graph with a marking is called marked graph. To denote it we write $L(G)=\{ x_1,...,x_n \}$. A morphism of marked graphs  $\phi:G\rightarrow G'$ is a morphism of graphs which preserves the marking, i.e. $\phi(x_i)=x_i'$ for every $i$ from 1 to $n$.

 We recall the notion of input datum given in \cite{HASSETT2003316}.
\begin{mydef}
Let $n\geq 1$. A weight datum is a $n$-tuple $\mathcal{A}=(a_1,...,a_n)\in ((0,1]\cap \mathbb{Q})^n.$ 
Given $g\geq 0$ an integer, an input datum is a pair $(g,\mathcal{A})$ with $\mathcal{A}$ being a weight datum such that $2g-2+a_1+...+a_n>0$. The integer $n$ is also called the length of the weight datum.
\end{mydef}
 The domain of all admissible weight data for genus $g$ and length $n$ is $$\mathcal{D}_{g,n}:=\lbrace (a_1,...,a_n)\in ((0,1]\cap \mathbb{Q})^n \text{ such that } a_1+...+a_n>2-2g\rbrace.$$
Note that for a fixed $n$ this space is $\mathcal{D}_{g,n}=((0,1]\cap \mathbb{Q})^n$ for every $g\geq 1$. We call $\mathcal{D}_{g,n}$ the space of stability conditions. We write $\overline{\mathcal{D}}_{g,n}$ for the space including also real coefficients, i.e. 
$$\overline{\mathcal{D}}_{g,n}=(0,1]^n,$$
as we need it for technical purposes.
 Let $G$ be a weighted marked graph, and let $L(G)=\lbrace x_1,...,x_n \rbrace$ be the set of its legs. For every $v\in V(G)$, we denote with $L(v)$ the set of legs adjacent to $v$, and with $\vert v\vert_E$ the number of half-edges that are not legs adjacent to $v$, i.e. $$ \vert v\vert_E=\text{val}(v)-\vert L(v)\vert .$$
For a given weight datum $\mathcal{A}$, we set also $$\vert v\vert_{\mathcal{A}}=\sum\limits_{x_i\in L(v)}a_i.$$
\begin{mydef}(\cite{ulirsch14tropical}, Definition 2.1)
Let $(g,\mathcal{A})$ be an input datum, $\mathcal{A}=(a_1,...,a_n)$, and $G$ a weighted marked graph with $n$ legs. We say that $G$ is stable of type $(g,\mathcal{A})$ (or that it is $(g,\mathcal{A})$-stable) if it is of genus $g$ and if for every vertex $v\in V(G)$
$$2w(v)-2+\vert v\vert_E+\vert v\vert_{\mathcal{A}}>0.$$
\end{mydef}
When $\mathcal{A}=(1,...,1):=1^{(n)}$, we recover the usual notion of stability for marked graphs, as the sum $|v|_E+|v|_{\mathcal{A}}$ becomes equal to $val(v)$.
\begin{myrem}
A weighted contraction of a stable graph of type $(g,\mathcal{A})$ is still a $(g,\mathcal{A})$-stable graph. 
\end{myrem}

\subsection{Tropical curves}
 \begin{mydef}
 An $n$-marked tropical curve of genus $g$ is a pair $\Gamma:=(G, l)$ where $G$ is a weighted marked graph of genus $g$ with $n$ legs and $l$ is a function $$l: E(G)\cup L(G)\rightarrow \mathbb{R}_{>0}\cup \lbrace \infty \rbrace $$ such that $l(x) = \infty$
if and only if $x$ is a leg or an edge adjacent to a vertex of valence 1 and weight 0.
\end{mydef}

 Let $w$ be the vertex weight function of $G$; if $w(v) = 0$ for every $v \in V (G)$, we write $w = \overline{0}$ and we say that the tropical curve is pure. A tropical curve is called regular if it is pure and if $G$ is a 3-regular graph, i.e. all of its vertices have valence 3. 

 The volume of a tropical curve is defined as the sum of its edge lenghts.
 
 We usually refer to $G$ as the underlying graph of the tropical curve, and we denote it with $G(\Gamma)$ when it is necessary. Legs are called marked points of the tropical curve. We also write $V(\Gamma)$, $E(\Gamma)$ and so on to indicate vertices, edges and other characteristics of the tropical curve, meaning the ones of the underlying graph. We also forget to recall the adjective $n$-marked when it is clear by the context. 
\begin{mydef}
A tropical curve is $(g,\mathcal{A})$-stable if its underlying graph is $(g,\mathcal{A})$-stable.
\end{mydef}

 We say that two tropical curves $\Gamma$ and $\Gamma'$ with $n$ marked points $L(\Gamma) = \lbrace x_1 , . . . , x_n \rbrace$ and $L(\Gamma ') = \lbrace x'_1 , . . . , x'_n \rbrace$ are isomorphic if there exists
an isomorphism of weighted marked graphs $\alpha$ between the underlying graphs $G(\Gamma)$ and $G(\Gamma ')$ such that $l(e) = l'(\alpha(e))$. We denote by $Aut(\Gamma)$ the group of automorphisms of a tropical curve $\Gamma$. Note that $Aut(\Gamma)\subset Aut(G(\Gamma))$.
\begin{mydef}
An $n$-marked extended tropical curve of genus $g$ is a pair $\Gamma:=(G, l)$ where $G$ is a weighted marked graph of genus $g$ and $l$ is a function $$l: E(G)\cup L(G)\rightarrow \mathbb{R}_{>0}\cup \lbrace \infty \rbrace $$ such that $l(x) = \infty$
for every leg.\newline
\end{mydef}
 Notice that the only difference with the not-extended tropical curves is that we possibly have infinite length edges. All the other notions, as e.g. the notion of tropical equivalence, remain unaltered. Extended tropical curves are needed to compactify the moduli space of tropical curves.
\subsection{Moduli Spaces}\label{moduli} For every $\mathcal{A}\in \mathcal{D}_{g,n}$, there are moduli spaces of tropical curves $M^{trop}_{g,\mathcal{A}}$ and extended tropical curves $\overline{M}^{trop}_{g,\mathcal{A}}$ which carry respectively the structure of a generalized
cone complex and of generalized extended cone complex, both in the sense of \cite{ACP}. Their construction is due to Ulirsch in \cite{ulirsch14tropical} generalizing previous constructions  of moduli spaces of tropical curves of \cite{BMV} and \cite{article}. Moreover, it is possible to define the locus $\Delta_{g,\mathcal{A}}$ of  $(g,\mathcal{A})$-stable tropical curves of genus $g$ with volume 1. 

 Consider the category $\mathcal{G}_{g,\mathcal{A}}$ of isomorphism classes of $(g,\mathcal{A})$-stable marked decorated graphs with morphisms generated by isomorphisms and contractions.
We define a natural contravariant functor:
$$\Sigma : \mathcal{G}_{g,\mathcal{A}}\rightarrow \mathbf{RPCC}$$
on the category $\mathbf{RPCC}$ of rational polyhedral cone complexes (see \cite{ulirsch14tropical} for further details) as follows: to each isomorphism class of $(g,\mathcal{A})$-stable marked decorated graph $G$ we associate the rational polyhedral cone $\sigma_G = \mathbb{R}_{\geq 0}^{|E(G)|}$. A weighted edge contraction $\pi: G \rightarrow G'$ induces the natural embedding $i_{\pi}: \sigma_{G'}\rightarrow \sigma_G$ of a face of $\sigma_G$. An automorphism of $G$ induces an automorphism of $\sigma_G$ if it is trivial, otherwise it induces a self-gluing.
Similarly there is also a natural functor $\overline{\Sigma}$ from $\mathcal{G}_{g,\mathcal{A}}$ into the category of extended rational polyhedral cone complexes that is given by sending $G$ into $\overline{\sigma}_G = \overline{\mathbb{R}}_{\geq 0}^{|E(G)|}$, where $\overline{\mathbb{R}}_{\geq 0}$ is $\mathbb{R}_{\geq 0}\cup \{\infty\}$.  The moduli space $M^{trop}_{g,\mathcal{A}}$ of $(g,\mathcal{A})$-stable tropical curves  is defined to be the colimit
$$M^{trop}_{g,\mathcal{A}}:=\lim_{\to} \sigma_G$$
taken over ($\mathcal{G}_{g,\mathcal{A}})^{op}$. The moduli space $\overline{M}^{trop}_{g,\mathcal{A}}$ of $(g,\mathcal{A})$-stable extended tropical curves is defined analogously using the compactified cones $\overline{\sigma}_G$'s. The moduli space $\Delta_{g,\mathcal{A}}$ is obviously a subspace of $M_{g,\mathcal{A}}^{trop}$, but it can be defined in the same way we did for the other spaces. Let $$\pi_G:= \left\{\ell: E(\G) \to \R_{\geq 0} \mid \sum_{e \in E(\G)} \ell(e) = 1 \right\} \subset \R^{E(G)}_{\geq 0}.$$
Then we have $$\Delta_{g,\mathcal{A}}:=\lim_{\to} \pi_G,$$ where the limit is taken again over ($\mathcal{G}_{g,\mathcal{A}})^{op}$.

\subsection{Graph Complexes}\label{complexes} Let $(g,\mathcal{A})$ be an input datum. The chain complex $G^{(g,\mathcal{A})}$ is a complex of rational vector spaces generated by elements $[\G,\omega]$ where $\G$ is an $n$-marked pure $(g,\mathcal{A})$-stable graph and $\omega$ is a total order on the set of the edges of $\G$. Generators are subject to the relation $$[\G,\omega]=sgn(\sigma)[\G',\omega']$$
if there is an isomorphism of $n$-marked graphs $\G\cong\G'$ under which the orders $\omega$ and $\omega'$ are related by the permutation $\sigma\in S_{|E(\G)|}$. This forces $[\G,\omega] = 0$ when $\G$ admits an automorphism that induces an odd permutation on the edges. The homological degree of $\G$ is $|E(\G)|-2g$.

So $G^{(g,\mathcal{A})}=\bigoplus G^{(g,\mathcal{A})}_j$, where
$$G^{(g,\mathcal{A})}_j=\{\text{Rational vector space spanned by elements }[\G,\omega] \text{ were } \G \text{ has }j+2g\text{ edges.}\}$$

The differential on $[\G,\omega] \neq 0$ is defined as $$ \partial [\G,\omega]=\sum_{i=0}^{|E(\G)|}(-1)^i[\G/e_i,\omega|_{E(\G)\setminus \lbrace e_i\rbrace}],$$
where $\G/e_i$ indicates the contraction and $\omega|_{E(\G)\setminus \lbrace e_i\rbrace}$ is the induced order. If $e_i$ is a loop, we interpret the corresponding term in the formula of the differential as zero.
\begin{myrem}
This generalizes the notions and the theory developed in \cite{CGP1} and \cite{CGP2}. In fact, when $\mathcal{A}=1^{(n)}$, we recover the definition of $G^{(g,n)}$ we have in \cite{CGP2}.
\end{myrem}

\section{Wall crossing for weight data}\label{section3}
 In this section we study the stability conditions on their own to show properties which reflect some phenomena of tropical curves depending only on the markings and the stability type.

 A chamber decomposition of $\mathcal{D}_{g,n}$ consists of a finite set $\mathit{W}$ of hyperplanes, called walls of the decomposition.
The chambers of the decomposition are the connected components of 
$$\mathcal{D}_{g,n}-\bigcup_{w\in\mathit{W}}w.$$

We will consider the fine chamber decomposition defined in \cite{HASSETT2003316}:
$$\mathit{W}_f=\left\lbrace \sum_{j\in S}a_j=1:S\subset \lbrace 1,...,n \rbrace, 2\leq|S|\leq n-2\delta_{g,0} \right\rbrace$$
where $\delta_{i,j}$ is the the Kronecker delta, and we will denote the set of the chambers with $\mathbf{K}$ (note that $W_f$ and $\mathbf{K}$ both depend on $g$ and $n$, but we avoid to recall this to lighten the notation).

Each wall in $\mathit{W}_f$ divides $\mathcal{D}_{g,n}$ in two connected components defined by inequalities \begin{center}
$a_{i_1}+...+a_{i_m}>1$ and $a_{i_1}+...+a_{i_m}<1$
\end{center}
 So each chamber $Ch$ is defined by a set of inequalities: each inequality is associated to a wall, and indicates in which of the two components determined by that wall the elements of $Ch$ lie.
\begin{myex}\label{Dg2}
Let $g\geq 1$ and $n=2$. We have $$\mathcal{D}_{g,2}=\lbrace (a_1,a_2)\in((0,1]\cup\mathbb{Q})^2 : 0<a_i\leq 1 \rbrace\subset(0,1]^2.$$
 Since $n=2$, we have only a possible $S$, i.e. only a wall given by $a_1+a_2=1$, and we get only two chambers as shown in Figure \ref{fig:Dg2}.
\begin{center}
\begin{figure}[h]
\scalebox{0.75}{
\begin{tikzpicture}
\draw[dashed, black, -latex](0,0)--(9,0);
node[below left];

\node at (8.2,0.3){1};
\node at (-0.3,8){1};
\draw[dashed, black, -latex](0,0)--(0,9);
 node[below left];
 \draw[thick, black] (8,0)--(8,8);
 \draw[thick, black] (0,8)--(8,8);
 \draw[thick, black] (8,0)--(0,8);
 \node at (5.8,4){$w=\lbrace a_1+a_2=1\rbrace$};
 \node at (3,2){ $a_1+a_2< 1$};
 \node at (5,6){$a_1+a_2>1$};
\end{tikzpicture}}
    \caption{The space $\mathcal{D}_{g,2}$ for $g\geq 1$.}
    \label{fig:Dg2}
\end{figure}
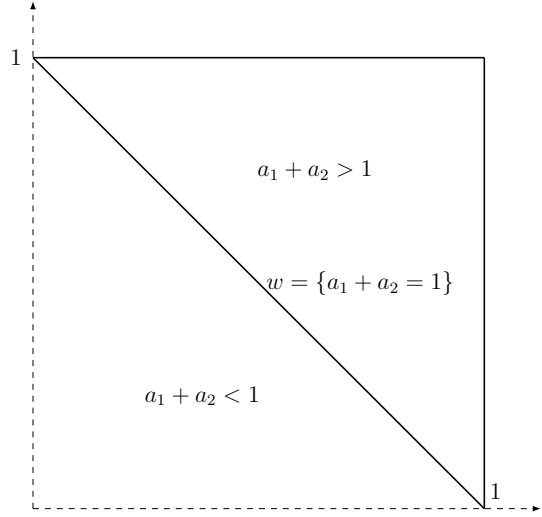
\end{center}
When $g=0$, we have the condition $a_1+a_2>2$, which is impossible since the $a_i$ are smaller or equal than 1, so there are not admissible input data. This  reflects the fact that there are no stable rational curves (either tropical or algebraic) with only two marked points.
\end{myex}
\begin{myex}\label{Dg3}
Let $g\geq 1$ and $n=3$. Here $\mathcal{D}_{g,3}\subset(0,1]^3$. We have $$W_f=\lbrace \lbrace a_1+a_2=1\rbrace , \lbrace a_1+a_3=1\rbrace , \lbrace a_2+a_3=1\rbrace , \lbrace a_1+a_2+a_3=1\rbrace \rbrace,$$ which are all planes in $\mathbb{R}^3$. The chambers of the fine decomposition are defined by the following sets of inequalities:

$\\ 
Ch_1:=\begin{cases}
a_1+a_2>1\\
a_1+a_3>1\\
a_2+a_3>1\\
a_1+a_2+a_3>1 \text{(which is implied by the three above)}
\end{cases}  \\ \\ \\
Ch_2:=\begin{cases}
a_1+a_2< 1\\
a_1+a_3>1\\
a_2+a_3>1\\
a_1+a_2+a_3>1 
\end{cases}; 
Ch_3:=\begin{cases}
a_1+a_2>1\\
a_1+a_3< 1\\
a_2+a_3>1\\
a_1+a_2+a_3>1 
\end{cases} \\ \\ \\
Ch_4:=\begin{cases}
a_1+a_2>1\\
a_1+a_3>1\\
a_2+a_3< 1\\
a_1+a_2+a_3>1 
\end{cases} ; 
Ch_5:=\begin{cases}
a_1+a_2< 1\\
a_1+a_3< 1\\
a_2+a_3>1\\
a_1+a_2+a_3>1 
\end{cases} \\ \\ \\
Ch_6:=\begin{cases}
a_1+a_2<1\\
a_1+a_3>1\\
a_2+a_3< 1\\
a_1+a_2+a_3>1 
\end{cases} ; 
Ch_7:=\begin{cases}
a_1+a_2>1\\
a_1+a_3< 1\\
a_2+a_3< 1\\
a_1+a_2+a_3>1 
\end{cases}\\ \\ \\
Ch_8:=\begin{cases}
a_1+a_2< 1\\
a_1+a_3< 1\\
a_2+a_3< 1\\
a_1+a_2+a_3>1 
\end{cases};
Ch_9:=\begin{cases}
a_1+a_2+a_3< 1 \text{(which implies the three below)}\\
a_1+a_2< 1\\
a_1+a_3< 1\\
a_2+a_3< 1\\
\end{cases}.\\ \\ \\ $
When $g=0$, we have $$\mathcal{D}_{0,3}=\lbrace (a_1,a_2,a_3)\in \mathbb{R}^3:0<a_i\leq 1, a_1+a_2+a_3>2\rbrace,$$ and  $\mathit{W}_f$ becomes the empty set, since the condition $2\leq |S|\leq 1$ is impossible. So there is only a non-empty chamber without walls.
\end{myex}
 From Example \ref{Dg3}, we can notice that there are inequalities which are not independent. Namely if $\sum_{j\in S}a_j<1$ then it has to be that $\sum_{j\in S'}a_j<1$ for every $S'\subset S$. At the same time if there is $S'\subset S$ such that $\sum_{j\in S'}a_j>1$, then $\sum_{j\in S}a_j>1$.
Moreover, whenever we have two sets $S$ and $T$ such that $S\cap T=\emptyset$ and both $\sum_{j\in S}a_j>1$ and $\sum_{j\in T}a_j>1$, then clearly $\sum_{j\in S\cup T}a_j>2$, so if there is a $T'\subset S\cup T$ such that $\sum_{j\in T'}a_j<1$, then $\sum_{j\in (S\cup T)-T'}a_j>1$.
\begin{myex}
Suppose we have $n=4$ and the defining inequalities $a_1+a_3>1$, $a_2+a_4>1$, and $a_1+a_2<1$. Then $a_1+a_2+a_3+a_4>2$, and so it must be that $a_3+a_4>1$, i.e. this inequality is forced by the others, otherwise we obtain an empty chamber.
\end{myex}

\subsection*{Chambers up to symmetry} Let us now consider the action of $S_n$ on $\mathcal{D}_{g,n}$ given by permuting the entries of weight data: 
\begin{center}
\begin{tikzcd}[row sep=tiny]
S_n\times \mathcal{D}_{g,n}\arrow[r] &\mathcal{D}_{g,n}\\
(\sigma, (a_1,...,a_n))\arrow[r, mapsto] & (a_{\sigma(1)},...,a_{\sigma(n)}).
\end{tikzcd}
\end{center}
\begin{myex}\label{exSn}
Consider $g=1$, $n=3$, let $\mathcal{A}=(\varepsilon, \frac{2}{3}, \frac{2}{3})\in \mathcal{D}_{1,3}$, for some $\epsilon>0$. Then the  graph $G$ in Figure \ref{fig:graphone} is $(1,\mathcal{A})$-stable: 
\begin{center}
\begin{figure}[h]
\scalebox{0.85}{
\begin{tikzpicture}
\draw[thick, black](0,0).. controls(-0.5,0.5) and (-0.5,-0.5) .. (0,0);
\draw[thick, black](0,0)--(1,0);
\draw[thick, black](0,0)--(0,1);
\node at (0.25,1){1};
\fill[black] (0,0) circle (0.4ex);
\fill[black] (1,0) circle (0.4ex);
\draw[thick, black](1,0)--(1.5,0.5);
\draw[thick, black](1,0)--(1.5,-0.5);
\node at (1.75,-0.5){3}; 
\node at (1.75,0.5){2};
\node at (0.5,0.25){$G$};
\end{tikzpicture}
}    
    \caption{A $(1,(\varepsilon,\frac{2}{3},\frac{2}{3}))$-stable graph.}
    \label{fig:graphone}
\end{figure}
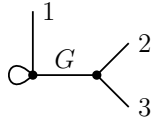
\end{center}

Let now $\sigma=(1\quad 2 \quad 3)\in S_3$, then $\sigma(\mathcal{A})=(\frac{2}{3},\varepsilon, \frac{2}{3})$, and one can observe that $G$ is not anymore stable with respect to the new weight datum. But changing the label of the legs of $G$ according to the same permutation $\sigma$ gives the new graph $G'$ in Figure \ref{fig:graphtwo}, which is $(1,\sigma(\mathcal{A}))$-stable.
\begin{center}
\begin{figure}[h]
\scalebox{0.85}{
\begin{tikzpicture}
\draw[thick, black](0,0).. controls(-0.5,0.5) and (-0.5,-0.5) .. (0,0);
\draw[thick, black](0,0)--(1,0);
\draw[thick, black](0,0)--(0,1);
\node at (0.25,1){2};
\fill[black] (0,0) circle (0.4ex);
\fill[black] (1,0) circle (0.4ex);
\draw[thick, black](1,0)--(1.5,0.5);
\draw[thick, black](1,0)--(1.5,-0.5);
\node at (1.75,-0.5){1}; 
\node at (1.75,0.5){3};
\node at (0.5,0.25){$G'$};
\end{tikzpicture}
}
    \caption{A $(1,(\frac{2}{3},\varepsilon,\frac{2}{3}))$-stable graph.}
    \label{fig:graphtwo}
\end{figure}
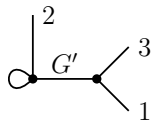
\end{center}
\end{myex}
The action of $S_n$ on $D_{g,n}$ descends on an action on $\mathbf{K}$, since a chamber is sent into another chamber by permutation of the coordinates. Under this action, two chambers are in the same orbit if there is a permutation in $S_n$ which sends all the inequalities defining the first chamber in the inequalities defining the second one by permuting the indices of the variables. We call the orbits of this action chambers up to symmetry. In particular, note that for a given $\sigma\in S_n$, $\sigma(Ch_{\mathcal{A}})=Ch_{\sigma(\mathcal{A})}$.

 We denote by $[Ch]$ the chamber up to symmetry of the chamber $Ch$, , and we denote the set of all the chambers up to symmetry by $[\mathbf{K}]$.
\begin{myex}\label{Dg3uts}
Let $g\geq1$, $n=3$. There are five orbits, namely
$$[\mathbf{K}]=\lbrace Ch_1\rbrace, \lbrace Ch_2, Ch_3, Ch_4\rbrace, \lbrace Ch_5, Ch_6, Ch_7\rbrace,\lbrace Ch_8\rbrace, \lbrace Ch_9\rbrace.$$
To see how to get an orbit, let us compute for example $\lbrace Ch_2, Ch_3, Ch_4\rbrace$. By the inequalities of Example \ref{Dg3}, we can see that $\sigma(Ch_2)=Ch_3$ for $\sigma=(2\quad 3),(1\quad 3\quad 2)$ and
$\tau(Ch_2)=Ch_4$ for $\tau=(1\quad 3),(1\quad 2\quad 3)$, while $id$ and $(1\quad 2)$ fix it. Analogously, we can see that $\theta(Ch_3)=Ch_4$ for $\theta=(1\quad 2), (1\quad 3\quad 2)$. To find a permutation $\sigma'$ such that $\sigma'(Ch_3)=Ch_2$ one can just pick the inverse of a permutation of one of the $sigma$'s, and analogously for $\tau'$ and $\theta'$.
\end{myex}

\begin{mydef}
Let $\mathcal{A}=(a_1,...,a_n)$ and $\mathcal{B}=(b_1,...,b_n)\in \mathcal{D}_{g,n}$ be two weight data such that $a_i\leq b_i$ for every $i$ from 1 to $n$. Then we write $\mathcal{A}\leq \mathcal{B}$.
\end{mydef} 
\begin{myrem}\label{rem:relation}
This relation is defined to reflect the following property. Let $(g,\mathcal{A}=(a_1,...a_n))$ and $(g,\mathcal{B}=(b_1,...b_n))$ be two input data such that $\mathcal{A}\leq \mathcal{B}$. Then a $(g,\mathcal{A})$-stable graph is always $(g,\mathcal{B})$-stable, because for every vertex $v\in V(G)$ we have $$0<2w(v)-2+|v|_E+|v|_{\mathcal{A}}\leq 2w(v)-2+|v|_E+|v|_{\mathcal{B}}.$$
In particular, for any weight data $\mathcal{A}$, a graph that is $(g,\mathcal{A})$-stable is always stable in the standard sense, since $\mathcal{A}\leq 1^{(n)}$ for every weight datum $\mathcal{A}$.
\end{myrem}

 We now focus on some technical results concerning weight data and chambers. We say that two chambers are adjacent if there is only a wall dividing them: in terms of inequalities, this means that there is only a subset $S\subset \lbrace 1,...,n \rbrace$ such that weight data in the first chamber satisfy $\sum_{i\in S}a_i<1$, and weight data in the second chamber satisfy the opposite inequality $\sum_{i\in S}a_i>1$, while for every other $S'\neq S$ the defining inequality it gives has the same direction for the two chambers or, in other words, they lie in the same half space induced by $S'$ on $\mathcal{D}_{g,n}$. We denote by $Ch_1|Ch_2$ the portion of the wall $w_S:=\{ \sum_{i\in S}a_i=1 \}\subset \overline{\mathcal{D}}_{g,n}$ dividing these two chambers, i.e. the subset of $w_S$ which verifies all the common defining inequalities of the two chambers.

\begin{mylemma} Let $Ch_1$ and $Ch_2$ be two adjacent chambers. The set $Ch1|Ch2$ is not empty.
\end{mylemma}

\begin{proof} 
Let $S$ be the subset of $\{1,...,n\}$ indexing the variables appearing in the defining inequalities which differ between $Ch_1$ and $Ch_2$. Suppose that weight data in $Ch_1$ verify $\sum_{i\in S}a_i<1$ and weight data in $Ch_2$ verify $\sum_{i\in S}a_i>1$. Let $\mathcal{X}=(x_1,...,x_n)$ in $Ch_1$ and $\mathcal{Y}=(y_1,...,y_n)$ in $Ch_2$. Consider the segment between $\mathcal{X}$ and $\mathcal{Y}$: each point of the segment can be described by
$$P_t=(1-t)\mathcal{X}+t\mathcal{Y}$$
for $t\in[0,1]$.

Consider the function $f:[0,1]\rightarrow \mathbb{R}$ sending $t$ into $(1-t)\sum_{i\in S}x_i+t\sum_{i\in S}y_i-1$: we have $f(0)=\sum_{i\in S}x_i-1<0$ since $\mathcal{X}$ is in $Ch_1$, while $f(1)=\sum_{i\in S}y_i-1>0$ since $\mathcal{Y}=(y_1,...,y_n)$ is in $Ch_2$. Then there must be a $t_0\in[0,1]$ such that $f(0)=0$, and this implies that $P_{t_0}$ belongs to the wall $w_S$. 

Now let $T\neq S$ be a subset of $\{1,...,n\}$, and suppose $\sum_{i\in T}a_i>1$ for each point in $Ch_1$ and $Ch_2$. Then $$(1-t)\sum_{i\in T}x_i+t\sum_{i\in T}y_i>(1-t)+t=1,$$ for every $t$ from 0 to 1.

Analogously if we pick $T'\neq S$ such that $\sum_{i\in T'}a_i<1$ for each point in $Ch_1$ and $Ch_2$, then 
$$(1-t)\sum_{i\in T}x_i+t\sum_{i\in T}y_i<(1-t)+t=1,$$ again for every $t$ from 0 to 1.

So in particular $P_{t_0}$ verifies all the common defining inequalities of the two chambers and belongs to the wall $w_S$, hence it lies in $Ch_1|Ch_2$.
\end{proof} 
\begin{myprop}\label{myprop:geqProp}
Let $Ch_1$ and $Ch_2$ be two adjacent chambers with  $S\subset \lbrace 1,...,n \rbrace$ such that data in $Ch_1$ satisfy $\sum_{i\in S}a_i<1$, while data in $Ch_2$ satisfy $\sum_{i\in S}a_i>1$, and for every other $T\neq S$ the corresponding inequalities agree for both chambers. Then there are $\mathcal{A}\in Ch_1$ and $\mathcal{B}\in Ch_2$ such that $\mathcal{A}\leq \mathcal{B}$. 
\end{myprop} 

\begin{proof}
By the previous Lemma the set $Ch_1|Ch_2$ is non empty, so we can choose $\mathcal{X}:=(x_1,...,x_n)\in Ch_1|Ch_2$. Let $\varepsilon$ being a number which is strictly less than $(\min\limits_{T\neq S}|\sum_{i\in S'}x_i-1|)$. Then we can pick $\mathcal{A}=(x_1,...,x_i-\varepsilon,...,x_n)$ and $\mathcal{B}=(x_1,...,x_i+\varepsilon,...,x_n)$ for some $i\in S$. 

Assume both have all rational components, then these two weight data verify all the common defining inequalities of $Ch_1$ and $Ch_2$. Indeed, if we choose a subset $T\neq S$ of $\{1,...,n\}$ such that $\sum_{i\in T}x_i<1$, then $\sum_{i\in T}a_i=\sum_{i\in S}x_i-\varepsilon<1$
and $\sum_{i\in T}b_i=\sum_{i\in T}x_i+\varepsilon<1$ since $\sum_{i\in T}x_i<1-\varepsilon$ by how we pick $\varepsilon$.

 Analogously if we choose $T\neq S$ such that $\sum_{i\in T}x_i>1$, then $\sum_{i\in T}b_i=\sum_{i\in T}x_i+\varepsilon>1$
and $\sum_{i\in S}a_i=\sum_{i\in T}x_i-\varepsilon>1$ since $\sum_{i\in T}x_i>1+\varepsilon$, so $\mathcal{A}$ and $\mathcal{B}$ verify all the common defining inequalities of $Ch_1$ and $Ch_2$.
Moreover $\sum_{i\in S}a_i<1$ and $\sum_{i\in S}b_i>1$, so it follows that $\mathcal{A}\in Ch_1,$ $\mathcal{B}\in Ch_2$, and $\mathcal{A}\leq \mathcal{B}$ by construction. 

If $\mathcal{A}=(x_1,...,x_i-\varepsilon,...,x_n)$ has irrational components, we can find a weight datum $\mathcal{A}'$ in $Ch_1$ approximating down each the irrational components $x_i$ with a rational number $x_i'$ such that $x_i-x_i'<\varepsilon$. The same can be done with $\mathcal{B}$ approximating up, so that we have $\mathcal{B}'\in Ch_2$ and by construction $\mathcal{A}'\leq \mathcal{B}'$.
\end{proof}

\begin{mydef}\label{orderonchambers}
Let $Ch_1, Ch_2\in \mathbf{K}$. We say that $Ch_1\leq Ch_2$ if they are equal or there are $S_1,...,S_t\subset \lbrace 1,...,n\rbrace$ such that for every $\mathcal{A}\in Ch_1$ we have $\sum_{i\in S_j}a_i<1$ and for every $\mathcal{B}\in Ch_2$ we have $\sum_{i\in S_j}b_i>1$ for every $j$ from 1 to $t$, and given any $S'\neq S_j$ for every $j$ from 1 to $t$, the two chambers are in the same half-space induced by the wall $\{\sum_{i\in S'}a_i=1\}$ on $\mathcal{D}_{g,n}$, i.e. all the other defining inequalities agree.
\end{mydef}   
\begin{myprop} 
The relation defined in \ref{orderonchambers} is a partial order on $\mathbf{K}$.
\end{myprop}\label{partialord}
\begin{proof}
The relation is reflexive by definition.\ It is also clearly transitive: let $Ch_1\leq Ch_2$ and $Ch_2 \leq Ch_3$, and let $S_1,...,S_t$ be the subsets of $\lbrace 1,...,n \rbrace $ on which the defining inequalities of $Ch_1$ and $Ch_3$ disagree. Fix a $j$ from $1$ to $t$. If the inequalities corresponding to $S_j$ agree for $Ch_1$ and $Ch_2$, then it has to be that $\sum_{i\in S_j}a_i<1$ for both, since it has to disagree with the inequality corresponding to $S_j$ of $Ch_3$ and $Ch_2\leq Ch_3$ by hypothesis. If the inequalities corresponding to $S_j$ disagree for $Ch_1$ and $Ch_2$, then it has to be that $\sum_{i\in S_j}a_i<1$ for $Ch_1$ and $\sum_{i\in S_j}a_i>1$ for $Ch_2$, since we have $Ch_1\leq Ch_2$. But then it has to be that $\sum_{i\in S_j}a_i>1$ also for $Ch_3$, otherwise it can not be that $Ch_2\leq Ch_3$, and so we have $Ch_1\leq Ch_3$.

For antisymmetry, let $Ch_1\leq Ch_2$ and $Ch_2\leq Ch_1$, and suppose they are different. By the first relation we get that there are $S_1,...,S_q$ such that their inequalities agree for $S'\neq S_j$, for $j=1,...,q$ and for them $\sum_{i\in S_j}a_i<1$ in $Ch_1$ and $\sum_{i\in S_j}a_i>1$ in $Ch_2$. On the other hand, the second equation gives that there are $T_1,...,T_r$ such that their inequalities agree for $S'\neq T_j$, and for them $\sum_{i\in T_j}a_i<1$ in $Ch_2$ and $\sum_{i\in T_j}a_i>1$ in $Ch_1$. But this is impossible, since such $T_j$'s can not exist by the first relation, so they must be the same chamber.
\end{proof}
The partial order we just defined on $\mathbf{K}$ naturally induces an order on $[\mathbf{K}]$.
\begin{mydef}\label{ord2}
Let $[Ch_1],[Ch_2]\in[\mathbf{K}]$ $[Ch_1]\leq[Ch_2]$ if there are two chambers $Ch_1\in[Ch_1]$ and $Ch_2\in[Ch_2]$ such that $Ch_1\leq Ch_2$.
\end{mydef}
 Then by Proposition \ref{partialord} we have the following.
\begin{mycor} The relation defined in \ref{ord2} is a partial order on the set $[\mathbf{K}]$. 
\end{mycor}

 Both the number of chambers and the number of chambers up to symmetry are finite, as showed in \cite{Enum}. For every $g$, there is a unique maximal chamber, namely the one given by the weight datum $1^{(n)}$. It is easy to see that it is the only element in the orbit of the action of $S_n$, and consequently $[Ch_{1^{(n)}}]$ is also the maximum with respect to the partial order on chambers up to symmetry.
When $g\geq 1$, there is also a minimal chamber (up to symmetry)  which contains all the admissible weight data $\mathcal{A}\leq\frac{1}{n}^{(n)}$. Also in this case, the orbit of the minimal chamber  is made only by itself. 

 From now on when $g\geq 1$, we denote by $\mathcal{E}$ a generic weight datum in the minimal chamber, which is denoted by $Ch_{\mathcal{E}}$.


\subsection*{An algorithm to compare weight data} 
Given two arbitrary weight data $\mathcal{A}=(a_1,...,a_n)$ and $\mathcal{B}=(b_1,...,b_n)$ in $\mathcal{D}_{g,n}$, it is possible to describe an algorithm which compares their chambers up to symmetry with respect to the order we put on them, and also says whenever they are not comparable. The procedure is the following:
\begin{itemize}
\item[1)] For a given input $n$, consider the group of permutations $S_n=\{\sigma_1=id,...,\sigma_k,...,\sigma_{n!}\}$. Consider also the inputs $\mathcal{A}$ and $\mathcal{B}$, and denote by $\mathcal{A}_k=(a_{1,k},...,a_{n,k})$ the datum $\sigma_k(\mathcal{A})$.
\item[2)] For each set $S\subset \{1,...,n\}$ such that $2\leq |S|\leq n-2\delta_{g,0}$, we compute the sums $\sum_{i\in S}b_i$.
\end{itemize}
The index $k$ will count the iterations of the algorithm. So to start the iteration here we set $k=1$ and $\mathcal{A}_1=\mathcal{A}$.

\begin{itemize}
    \item[3)] For each set $S\subset \{1,...,n\}$ such that $2\leq |S|\leq n-2\delta_{g,0}$, we compute the sums $\sum_{i\in S}a_{i,k}$.
\end{itemize}
Then we can have the following outputs:
\begin{itemize}
    \item[3.1)] If the condition $\sum_{i\in S}a_{i,k}\leq1$ holds if and only if the condition $\sum_{i\in S}b_i\leq 1$ holds (and of course  $\sum_{i\in S}a_{i,k}>1$ if and only if $\sum_{i\in S}b_i>1$), then $[Ch_\mathcal{A}]=[Ch_{\mathcal{B}}]$; 
    \item[3.2)] If there are $S_1,...,S_d$ such that $\sum_{i\in S_j}a_{i,k}\leq1$ and $\sum_{i\in S_j}b_i>1$, while for every other $S\neq S_j$ for every $j=1,...,d$ we have $\sum_{i\in S}a_{i,k}\leq1$ if and only if $\sum_{i\in S}b_i\leq 1$, then $[Ch_\mathcal{A}]\leq[Ch_{\mathcal{B}}]$. If the same happens but with the roles of $\mathcal{A}_k$ and $\mathcal{B}$ are reversed then $[Ch_\mathcal{B}]\leq[Ch_{\mathcal{A}}]$
    \item[3.3)]If there are $S$, $T$ such that $\sum_{i\in S}a_{i,k}\leq1$ and $\sum_{i\in T}b_i>1$ while $\sum_{i\in T}a_{i,k}>1$ and $\sum_{i\in T}b_i\leq1$, we let the index $k$ grow by one.
    
    If $k\leq n!$ we restart from the point $(3)$ of the algorithm. 
    
    When $k=n!+1$, we can conclude that $[Ch_\mathcal{A}]$ and $[Ch_{\mathcal{B}}]$ are not comparable in the partial order.
\end{itemize}
When $n$ grows, this algorithm is not really efficient as it needs an extremely large number of operations: in the worst case, if $g\geq 1$ we have to compute $(2^n-n)(n!+1)$ sums and make $(2^n-n)n!$ comparisons of results, while if $g=0$ the number of sums is $(2^n-2n-1)(n!+1)$, with $(2^n-2n-1)n!$ comparisons.

Notice that in the case of comparable chambers up to symmetry a smarter choice of the permutation can reduce the number of the operations, but they will never fall below the number of operations of the best case, which is when we need a single iteration. 
In this case the number of sums is $2(2^n-n)$  if $g\geq 1$ and $2(2^n-2n-1)$ if $g=0$, while the number of comparisons is $(2^n-n)$ if $g\geq 1$ and $(2^n-2n-1)$ if $g=0$.

\section{Proof of the first main theorem}\label{section4}
In this section we give the proof
of Theorem \ref{Main1}, restated here.
\begin{theorem*}
Let $g\geq 0$, $n\geq 1$ be two integers, and $\mathcal{A}\in \mathcal{D}_{g,n}$ a weight datum.
There are filtrations of $M^{trop}_{g,\mathcal{A}}$ given by embeddings induced by the partial order on the set of chambers up to simmetry. Namely, an ordered sequence  $$[Ch_{\mathcal{A}_1}]\leq [Ch_{\mathcal{A}_2}]\leq ...\leq [Ch_{\mathcal{A}_p}] \leq ... \leq [Ch_{\mathcal{A}_{N-1}}]\leq [Ch_{\mathcal{A}}]$$ induces a filtration is
    $$M^{trop}_{g,\mathcal{\mathcal{A}}_1}\hookrightarrow M^{trop}_{g,\mathcal{A}_2}\hookrightarrow  ...\hookrightarrow M^{trop}_{g,\mathcal{A}_p} \hookrightarrow ... \hookrightarrow M^{trop}_{g,\mathcal{A}_{N-1}}\hookrightarrow M^{trop}_{g,\mathcal{A}}. $$
 The same result holds if we replace $M_{g,\mathcal{A}}^{trop}$ with the moduli space of extended weighted tropical $(g,\mathcal{A})$-stable curves $\overline{M}_{g,\mathcal{A}}^{trop}$  or the moduli space of $(g,\mathcal{A})$-stable tropical curves with volume 1 $\Delta_{g,\mathcal{A}}$.
\end{theorem*}
The strategy of the proof is to show that there are filtrations of the graph categories $\mathcal{G}_{g,\mathcal{A}}$ induced by the partial order on the chambers up to symmetry. Then the direct limit description of the moduli spaces of tropical curves will give us the result.
\subsection{Wall-Crossing properties for graphs} Let $\mathcal{G}_{g,\mathcal{A}}$ be the graph category introduced in \cite{ulirsch14tropical}, where objects are $(g,\mathcal{A})$ weighted marked graphs with the morphisms described in Section \ref{section2}. We start by showing that weight data lying in the same chamber define equal categories of graphs.
\begin{mylemma}
\label{lemmaCoarse} Consider the map $$\Psi:\mathcal{D}_{g,n}\setminus\bigcup_{w\in \mathit{W}_f}w \rightarrow \{\mathcal{G}_{g,\mathcal{A}}|\mathcal{A} \text{ is a weight datum}\}$$ sending a weight datum $\mathcal{A}$ in the graph category $\mathcal{G}_{g,\mathcal{A}}$.
Then

1)The map $\Psi$ is constant on each chamber of the fine chamber decomposition.

 2)The fine chamber decomposition is the coarsest one with the above property, i.e. if $\mathcal{A}$ and $\mathcal{B}$ are in two different chambers of the fine chamber decomposition, their image under the above map is different.
\end{mylemma} 
\begin{proof} We follow the lines of Proposition 6.2 of \cite{ulirsch14tropical}, since the proof relies only on graphs. The map is clearly constant on the fine chambers of $\mathcal{D}_{g,n}$. It suffices to show that
$\mathcal{
G}_{g,\mathcal{A}}$ changes whenever we cross a wall. Let $S\subset\{1, . . . , n\}$ with $2\leq|S|\leq n$ be the subset indexing the equation of the wall $w=\{\sum_{i\in S}a_i=1\}$.

 Suppose first that $g\geq 1$. 
Let $S\subset\{1, . . . , n\}$ with $2\leq|S|\leq n$ be the subset indexing the equation of the wall $w=\{\sum_{i\in S}a_i=1\}$, and  consider the graph containing one edge between two vertices, one of weight 0 and one of weight $g$, with all the legs with index in $S$ being incident to the one with weight 0 (see Figure \ref{genusgraph}). Then the graph is stable of type $(g,\mathcal{A})$ if $\sum_{i\in S}a_i$ is greater than 1, otherwise it is not, i.e. changing the half-space of $\mathcal{D}_{g,n}$  we are also changing the category $\mathcal{G}_{g,\mathcal{A}}$, since the graphs are not the same.
\begin{center}
\begin{figure}[h]\scalebox{0.75}{
\begin{tikzpicture}
\draw[thick, black](0,0)--(3,0);
\draw[thick, black](0,0)--(1,0);
\draw[thick, black](0,0).. controls(-2,1.1) and (-2,-1.1) .. (0,0);
\draw[thick, black](0,0).. controls(-0.5,0.25) and (-0.5,-0.25) .. (0,0);
\draw[thick, black](0,0).. controls(-1,0.55) and (-1,-0.55) .. (0,0);
\draw[thick, black](0,0).. controls(-1.5,0.76) and (-1.5,-0.76) .. (0,0);
\draw (0,0) circle (0.4ex);
\fill[black] (0,0) circle (0.4ex);
\draw (3,0) circle (0.4ex);
\fill[black] (3,0) circle (0.4ex);
\node at (0,0.3){$u$};
\node at (3,0.3){$v$};
\draw[thick, black](3.5,1)--(3,0);
\draw[thick, black](4,0)--(3,0);
\draw[thick, black](3.5,-1)--(3,0);
\draw[dashed, black](3.75,0.5)--(3,0);
\draw[dashed, black](3.75,-0.5)--(3,0);
\node at (5.5,0.3){legs indexed by $S$};
\node at (-1,-0.6){$g$ loops};
\draw[thick, black](0,0)--(0.5,-0.7);
\draw[dashed, black](0,0)--(0.7,-0.5);
\draw[thick, black](0,0)--(0.5,0.7);
\draw[dashed, black](0,0)--(0.7,0.5);
\node at (1,-1){other legs};
\end{tikzpicture}
}
    \caption{}
    \label{genusgraph}
\end{figure}
\end{center}
 In the case $g = 0$, let $S\subset\{1, . . . , n\}$ with $2\leq|S|\leq n-2$ be the subset indexing the equation of the wall $w=\{\sum_{i\in S}a_i=1\}$, and consider the graph with two vertices of weight 0 connected by an edge and legs incident to the first vertex having indices in $S$, while the others are incident to the second vertex (Figure \ref{nogenus}).
\begin{center}
\begin{figure}[h]\scalebox{0.75}{
    
\begin{tikzpicture}
\draw[thick, black](0,0)--(3,0);
\draw (0,0) circle (0.4ex);
\fill[black] (0,0) circle (0.4ex);
\draw (3,0) circle (0.4ex);
\fill[black] (3,0) circle (0.4ex);
\node at (0,0.3){$u$};
\node at (3,0.3){$v$};
\draw[thick, black](3.5,1)--(3,0);
\draw[thick, black](4,0)--(3,0);
\draw[thick, black](3.5,-1)--(3,0);
\draw[dashed, black](3.75,0.5)--(3,0);
\draw[dashed, black](3.75,-0.5)--(3,0);
\node at (5.4,0.3){legs indexed by $S$};
\draw[thick, black](0,0)--(-0.5,1);
\draw[dashed, black](0,0)--(-0.75,-0.5);
\draw[dashed, black](0,0)--(-0.75,0.5);
\draw[thick, black](0,0)--(-0.5,-1);
\node at (1,-1){other legs};
\draw[thick, black](-1,0)--(0,0);
\end{tikzpicture}
}
\caption{}
    \label{nogenus}
\end{figure}
\end{center}
Suppose that $\sum_{i\in S} a_i\leq 1$, then the condition $\sum_{i=1}^n a_i > 2$ implies $\sum_{i\notin S} a_i > 1$. So when crossing the wall 
$\sum_{i\in S} a_i= 1$ without changing the $a_i$ such that $i\notin S$ we obtain that the described graph is stable of type $(0, \mathcal{A})$ if $\sum_{i\in S}a_i> 1$, and it is not otherwise, and again we conclude that changing the half-space of $\mathcal{D}_{g,n}$ the category $\mathcal{G}_{g,\mathcal{A}}$ is also different.
\end{proof}
 In the following lemma we will consider the case of weight data lying on walls.
\begin{mylemma}\label{wallprop}
Let $d\geq 1$ be an integer and $S_1,...,S_d\subset \{1,...,n\}$. Let $Ch$ be a chamber such that $\sum_{i\in S_j}a_i<1$ for every $j$ from 1 to $d$. Let $\mathcal{R}$ be a weight datum such that $\sum_{i\in S_j}r_i=1$, for every $j$ from 1 to $d$, while for every $S\neq S_j$ $\mathcal{R}$ belongs to the half-space induced by $S$ containing $Ch$. Then $\mathcal{G}_{g,\mathcal{A}}$ is equal to $\mathcal{G}_{g,\mathcal{R}}$, for every $\mathcal{A}\in Ch$.
\end{mylemma}
\begin{proof}
By definition, for all $S$ and for every $\mathcal{A}\in Ch$ we have $\sum_{i\in S}a_i<1$  if and only if $\sum_{i\in S}r_i\leq 1$. Indeed since $\mathcal{A}$ belongs to a Chamber $\sum_{i\in S}a_i\neq 1$ for any $S\subset\{1,...,n\}$. 
Moreover, we can find $\mathcal{A}\in Ch$ with the property that $\mathcal{A}\leq\mathcal{R}$. Indeed, let $\mathcal{S}=\bigcup_{j=1}^d S_i$ and let $\varepsilon:=(\min\limits_{S'\neq S}|\sum_{i\in S'}x_i-1|)/2|\mathcal{S}|$. We can consider $\mathcal{A}=(\overline{r}_1,...,\overline{r}_n)$ where:
\begin{center}
  $
\overline{r}_i:=\begin{cases}
r_i \text{ if }i\notin \mathcal{S}\\
r_i-\varepsilon \text{ if }i\in \mathcal{S}
\end{cases}
$  
\end{center}
 By construction if $\sum_{i\in S_j}r_i=1$ then $\sum_{i\in S_j}\overline{r}_i<1$ for every $j=1,...,d$. Moreover if $\sum_{i\in S}r_i<1$ then $\sum_{i\in S}\overline{r}_i<1$ and if $\sum_{i\in S}r_i>1$ then $\sum_{i\in S}\overline{r}_i>1$. 
 
By construction $\mathcal{A}\leq \mathcal{R}$, so the category $\mathcal{G}_{g,\mathcal{A}}$ is a full subcategory of $\mathcal{G}_{g,\mathcal{R}}$ by Remark \ref{rem:relation}. 
As in \ref{myprop:geqProp} we can assume $\mathcal{A}$ to have alla rational components, so that $\mathcal{A}\in Ch$

 Now by contradiction suppose there is $G \in Ob(\mathcal{G}_{g,\mathcal{R}} )\smallsetminus Ob(\mathcal{G}_{g,\mathcal{A}} )$. If this happens, there is $v\in V(G)$ such that $$2w(v)-2+|v|_E+|v|_{\mathcal{A}}<0<2w(v)-2+|v|_E+|v|_{\mathcal{R}}.$$
This implies that there is an $S$ such that $|v|_{\mathcal{A}}=\sum_{i\in S}a_i<1$ while $|v|_{\mathcal{R}}=\sum_{i\in S}r_i>1$, indeed $2w(v)-2+|v|_E+|v|_{\mathcal{A}}<0$ implies $2w(v)-2+|v|_E<-|v|_{\mathcal{A}}$, so $w(v)=0$ and $|v|_E\leq 1$, since they are all integers. But this is a contradiction by our hypothesis on $Ch$ and $\mathcal{R}$, so the result follows.
\end{proof}

\begin{myprop}
Let $Ch_1,Ch_2\in \mathbf{K}$ be two chambers, $\mathcal{A}\in Ch_1$, $\mathcal{B}\in Ch_2$ two weight data. Then, if $Ch_1\leq Ch_2$,
the category $\mathcal{G}_{g,\mathcal{A}}$ is a full subcategory of $\mathcal{G}_{g,\mathcal{B}}$.
\end{myprop}
\begin{proof} Suppose $Ch_1$ and $Ch_2$ are different, otherwise the result is trivial by Lemma \ref{lemmaCoarse}.
Since $Ch_1\leq Ch_2$, there are $S_1,..., S_d\subseteq \{1,...,n\}$ such that $\sum_{i\in S}a_i<1$ for $(a_1,...,a_n)\in Ch_1$ and $\sum_{i\in S}b_i>1$ for $(b_1,...,b_n)\in Ch_2$, for $j=1,...,d$, while for $S\neq S_j$ the defining inequalities of $Ch_1$ and $Ch_2$ have the same direction.

 Suppose first $d=1$. Since $\mathcal{G}_{g,\mathcal{A}}$'s are constant on each chamber, by Proposition \ref{myprop:geqProp} we can find $\mathcal{A}'\in Ch_1$, $\mathcal{B}'\in Ch_2$ such that $\mathcal{A}'\leq \mathcal{B}'$ and we have $\mathcal{G}_{g,\mathcal{A}}=\mathcal{G}_{g,\mathcal{A}'}$ and $\mathcal{G}_{g,\mathcal{B}}=\mathcal{G}_{g,\mathcal{B}'}$. The inclusion follows from Remark \ref{rem:relation}.

 Let $d\geq 2$, and suppose by induction that for every $d'<d$, given two chambers $Ch_1'$ and $Ch_2'$ which inequalities agree for all but $d'$ sets of indices $S_1',...,S_{d'}'$, and for every $\mathcal{A}\in Ch_1$ we have $\sum_{i\in S_j'}a_i<1$ while for every $\mathcal{B}\in Ch_2$ we have $\sum_{i\in S_j'}a_i>1$, then we can find a weight datum in $\mathcal{A'}\in Ch_1'$ and $\mathcal{B'}\in Ch_2'$ such that $\mathcal{A'}\leq\mathcal{B'}$.

 First, suppose there is a chamber $Ch_3$ which defining inequalities agree with the one of $Ch_1$ except for $S_1,...,S_c$, for a number $c<d$. Then by induction there is $\mathcal{A}\in Ch_1$, $\mathcal{P},\mathcal{Q}\in Ch_3$ and $\mathcal{B}\in Ch_2$ such that $\mathcal{A}\leq
\mathcal{P}$ and $\mathcal{Q}\leq \mathcal{B}$. Then we can induce inclusions $$\mathcal{G}_{g,\mathcal{A}}\hookrightarrow \mathcal{G}_{g,\mathcal{P}}=\mathcal{G}_{g,\mathcal{Q}}\hookrightarrow\mathcal{G}_{g,\mathcal{B}},$$ so by composition we get the desired inclusion.

 Suppose now there is no chamber $Ch_3$, i.e. for every weight datum $\mathcal{Q}$ not belonging to walls such that $\sum_{i\in S_j}q_i>1$ for some $i=1,...,d$, then $\sum_{i\in S_k}q_i>1$ for every $k=1,...,d$ and the same holds picking the symbol $<$ instead of the symbol $>$.

Let $\mathcal{X}\in Ch_1$, $\mathcal{Y}\in Ch_2$ and consider the segment which goes from $\mathcal{X}$ and $\mathcal{Y}$: then every point on the segment $P_t=(1-t)\mathcal{X}+t\mathcal{Y}$ for $t\in[0,1]$ verifies all the common inequalities, and there is at least a $t'\in [0,1]$ such that $\mathcal{P}_{t'}=(p_1,...,p_n)$ belongs to a wall $w_{S_k}$, for some $k=1,...,d$. Let $\varepsilon=min_{\{T|\sum_{i\in T}p_i\neq 1\}}|1-\sum_{i\in T}p_i|$ and let $\delta<\varepsilon$. Define $\mathcal{B}=(b_1,...,b_n)=(p_1+\frac{\delta}{n},...,p_n+\frac{\delta}{n})$. It follows by construction that $\mathcal{B}$ verifies all the inequalities verified by $\mathcal{P}_t$, and since $\sum_{i\in S_k}b_i>1$ and it can not be on a wall it belongs to $Ch_2$. Analogously we can define $\mathcal{A}=(a_1...,a_n)=(p_1-\frac{\delta}{n},...,p_n-\frac{\delta}{n})$. Since $\sum_{i\in S_k}a_i<1$ while all other equalities agree with the one of $\mathcal{P}_{t'}$ then it belongs to $Ch_1$, and $\mathcal{A}\leq \mathcal{B}$ by construction, so we conclude by observing that we can suppose them to have all rational components as in the proof of \ref{myprop:geqProp}.
\end{proof}
\begin{myrem}
The proof of this proposition also shows that for any two chambers $Ch_1\leq Ch_2$ we can find a weight datum in $\mathcal{A}\in Ch_1$ and $\mathcal{B}\in Ch_2$ such that $\mathcal{A}\leq\mathcal{B}$.
\end{myrem}
 By Example \ref{exSn} we can easily deduce that each time we have two weight data $\mathcal{A}$ and a permutation $\sigma \in S_n$, the two categories $\mathcal{G}_{g,\mathcal{A}}$ and $\mathcal{G}_{g,\sigma(\mathcal{A})}$ are isomorphic, since it is enough to send each graph to the one obtained relabeling its legs according to $\sigma$, without changing morphisms (they are equal if the chosen permutation acts trivially on the chamber). So, up to isomorphism, there is only a category $\mathcal{G}_{g,\mathcal{A}}$ for each chamber up to symmetry, which is the one containing the chamber to which $\mathcal{A}$ belongs. Therefore, the latter Proposition can be rephrased including this symmetry property, giving the following:
 \begin{myprop}
Let $[Ch_1]$ and $[Ch_2]$ be two chambers up to symmetry, $\mathcal{A}\in Ch_1\in [Ch_1]$, $\mathcal{B}\in Ch_2\in [Ch_2]$. There is an inclusion as a full subcategory $\mathcal{G}_{g,\mathcal{A}}\hookrightarrow \mathcal{G}_{g,\mathcal{B}}$ each time $[Ch_1]\leq [Ch_2]$. It is an isomorphism if the two chambers up to symmetry are the same.
\end{myprop}
 Since the stability conditions on tropical curves are defined on their underlying graphs, everything we showed so far can be easily generalized for $M_{g,\mathcal{A}}^{trop}$, $\overline{M}_{g,\mathcal{A}}^{trop}$ and $\Delta_{g,\mathcal{A}}$. We can resume everything in the following Proposition.

\begin{myprop}\label{conseq}
Let $g,n\geq 1$ be two integers and $\mathcal{A}$ and $\mathcal{B}$ two weight data in $\mathcal{D}_{g,n}$.
\begin{itemize}
    \item[1)]\label{prop1} If $\mathcal{A}\leq \mathcal{B}$, then $M_{g,\mathcal{A}}^{trop}\subset M_{g,\mathcal{B}}^{trop}$;
    \item[2)]\label{prop2} If $\mathcal{A}$ and $\mathcal{B}$ are in the same chamber, then $M_{g,\mathcal{A}}^{trop}$ is equal to $M_{g,\mathcal{B}}^{trop}$;
    \item[3)]\label{prop3} If $\mathcal{A}$ and $\mathcal{B}$ are obtained one from the other through a permutation of coordinates, then $M_{g,\mathcal{A}}^{trop}$ is homeomorphic to $M_{g,\mathcal{B}}^{trop}$ through a relabeling homeomorphism;
    \item[4)]\label{prop4}If $\mathcal{A}$ and $\mathcal{B}$ are in chambers which belong to the same orbit, then $M_{g,\mathcal{A}}^{trop}$ is homeomorphic to $M_{g,\mathcal{B}}^{trop}$ through a relabeling homeomorphism.
\end{itemize}
The same results hold if we replace the $M_{g,\mathcal{A}}^{trop}$'s with the $\overline{M}_{g,\mathcal{A}}^{trop}$'s or the $\Delta_{g,\mathcal{A}}$'s.
\end{myprop}
\begin{myrem}
A priori, given $\mathcal{A}$ and $\mathcal{B}$ in different chambers, we can not say if two moduli spaces $M_{g,\mathcal{A}}^{trop}$ and $M_{g,\mathcal{B}}^{trop}$ are not homeomorphic as topological spaces. The point of using the relabeling homeomorphism is that they preserve the moduli space structure.
\end{myrem}
\begin{proof}(Theorem \ref{Main1})
We can now conclude the proof of Theorem \ref{Main1}. Consider an ordered sequence of chambers up to symmetry $$[Ch_{\mathcal{A}_1}]\leq [Ch_{\mathcal{A}_2}]\leq ...\leq [Ch_{\mathcal{A}_p}] \leq ... \leq [Ch_{\mathcal{A}_{N-1}}]\leq [Ch_{\mathcal{A}}].$$ At each step of the filtration $[Ch_{\mathcal{A}_p}]\leq [Ch_{\mathcal{A}_{p+1}}]$ we can find two weight data $\mathcal{A}\leq \mathcal{B}$ and two chambers $Ch_1\leq Ch_2$ such that $\mathcal{A}\in Ch_1\in [Ch_{\mathcal{A}_p}]$ and $\mathcal{B}\in Ch_2\in [Ch_{\mathcal{A}_{p+1}}]$. Then $M_{g,\mathcal{A}_p}^{trop} \approx M_{g,\mathcal{A}}^{trop}\subset M_{g,\mathcal{B}}^{trop} \approx M_{g,\mathcal{A}_{p+1}}^{trop}$, where by $\approx$ we indicate homeomorphisms of the above proposition, \ref{prop3} and the inclusion is the one of \ref{prop1}. This induces the desired inclusion map $M_{g,\mathcal{A}_p}^{trop}\hookrightarrow M_{g,\mathcal{A}_{p+1}}^{trop}$.  The result of the Theorem then follows repeating this reasoning for each step of the sequence.
\end{proof}

\begin{myex}\label{final}
Consider the case $g\geq 1$, $n=3$.  We saw the chamber decomposition in Example \ref{Dg3}, and the chamber decomposition up to symmetry in Example \ref{Dg3uts}. We choose a weight datum for each chamber, and for symmetric chambers we choose data obtained after a permutation:
\begin{itemize}
\item $(1,1,1)\in Ch_1 ;$
\item $(\frac{12}{27},\frac{14}{27},1-\varepsilon)\in Ch_2 ;$
\item $(\frac{14}{27},1-\varepsilon, \frac{12}{27})\in Ch_3 ;$
\item $(1-\varepsilon, \frac{12}{27},\frac{14}{27})\in Ch_4; $
\item $(\frac{12}{27},\frac{14}{27}, \frac{14}{27}-\varepsilon)\in Ch_5,$ 
\item $(\frac{14}{27}-\varepsilon,\frac{12}{27},\frac{14}{27})\in Ch_6; $
\item $(\frac{14}{27}, \frac{14}{27}-\varepsilon, \frac{12}{27})\in Ch_7;$
\item $(\frac{4}{9}-\varepsilon,\frac{4}{9}-\varepsilon,\frac{4}{9}-\varepsilon)\in Ch_8;$
\item $(\frac{1}{3},\frac{1}{3},\frac{1}{3}-\varepsilon)\in Ch_9,$
\end{itemize}
for $0< \varepsilon <\frac{1}{27}$. We pick these $\varepsilon$ perturbations in order to have our weight data in the interior of chambers. Since moduli spaces are constant on each chamber, for every weight datum $\mathcal{A}\in\mathcal{D}_{g,3}$, the moduli space $M_{g,\mathcal{A}}^{trop}$ is the same as the space obtained picking one of the nine data above (the one which lies in the same chamber of $\mathcal{A}$). So by the partial order on the set of the chambers, we get the following diagram:
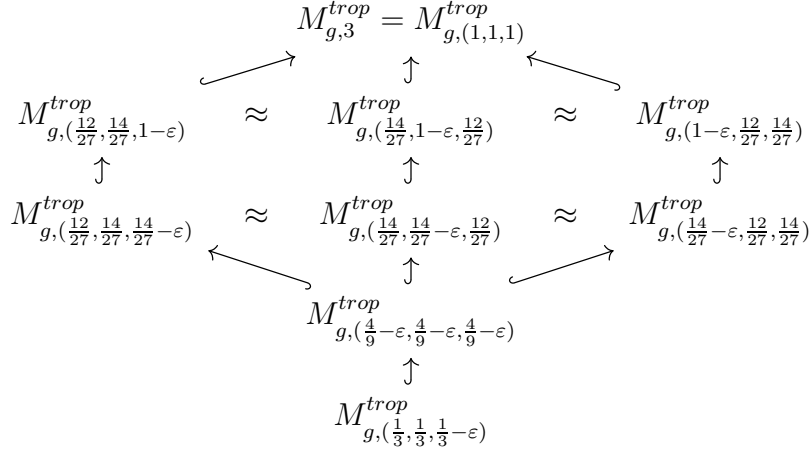
\begin{center}
\begin{tikzcd}[row sep=small]
& M_{g,3}^{trop}=M_{g,(1,1,1)}^{trop} \\
M_{g,(\frac{12}{27},\frac{14}{27},1-\varepsilon)}^{trop} \arrow[ur, hookrightarrow] \arrow[r, phantom, "\approx "] & M_{g,(\frac{14}{27},1-\varepsilon, \frac{12}{27})}^{trop} \arrow[u, hookrightarrow] \arrow[r, phantom, "\approx "] & M_{g,(1-\varepsilon, \frac{12}{27},\frac{14}{27})}^{trop} \arrow[ul, hookrightarrow] \\
M_{g,(\frac{12}{27},\frac{14}{27}, \frac{14}{27}-\varepsilon)}^{trop} \arrow[u, hookrightarrow] \arrow[r, phantom, "\approx "] & M_{g,(\frac{14}{27}, \frac{14}{27}-\varepsilon, \frac{12}{27})}^{trop}\arrow[u, hookrightarrow] \arrow[r, phantom, "\approx "] & M_{g,(\frac{14}{27}-\varepsilon,\frac{12}{27},\frac{14}{27})}^{trop} \arrow[u, hookrightarrow] \\
& M_{g,(\frac{4}{9}-\varepsilon,\frac{4}{9}-\varepsilon,\frac{4}{9}-\varepsilon)}^{trop}\arrow[u, hookrightarrow] \arrow[ur, hookrightarrow] \arrow[ul, hookrightarrow]  \\
& M_{g,(\frac{1}{3},\frac{1}{3},\frac{1}{3}-\varepsilon)}^{trop} \arrow[u, hookrightarrow]  
\end{tikzcd}
\end{center}
In the diagram, we indicated the relabeling homeomorphism of Proposition \ref{conseq} with $\approx$. 

 For example, the following filtration
$$ M_{g,(\frac{1}{3},\frac{1}{3},\frac{1}{3}-\varepsilon)}^{trop} \subset  M_{g,(\frac{4}{9}-\varepsilon,\frac{4}{9}-\varepsilon,\frac{4}{9}-\varepsilon)}^{trop} \subset M_{g,(\frac{14}{27}-\varepsilon, \frac{12}{27}, \frac{14}{27})}^{trop} \subset M_{g,(1-\varepsilon, \frac{12}{27},\frac{14}{27})}^{trop} \subset M_{g,3}^{trop},$$
of the space $M_{g,3}^{trop}$ can be obtained by picking the right column of the diagram.

 Clearly, the same diagram and the same filtrations also work for $\overline{M}_{g,\mathcal{A}}^{trop}$ and $\Delta_{g,\mathcal{A}}$. Notice that in this case the partial order induced on the chambers up to symmetry becomes total, since we can say for each couple of chambers up to symmetry which of them is greater or equal than the other. This is not true in general, as showed in the following example.
\end{myex}

\begin{myex}
Let $g\geq 0$, $n=8$. Consider the chamber $Ch_1$ defined by the following set of inequalities:
\begin{center}
  $
\begin{cases}
a_1+a_2>1\\
a_1+a_3>1\\
a_1+a_4>1\\
a_1+a_5>1\\
a_1+a_6>1\\
a_1+a_7>1\\
a_i+a_j<1 \text{ for any other couple of indices}\\
\sum_{i\in S}a_i>1 \text{ for any }S \text{ such that } |S|\geq 3
\end{cases}
$  
\end{center}
This is not empty, for example the datum $$\mathcal{A}_1=(\frac{1}{2}+2\varepsilon,\frac{1}{2}-\varepsilon,\frac{1}{2}-\varepsilon,\frac{1}{2}-\varepsilon,\frac{1}{2}-\varepsilon,\frac{1}{2}-\varepsilon,\frac{1}{2}-\varepsilon,2\varepsilon)$$ belongs to $Ch_1$ for a sufficiently small $\varepsilon$. Consider now the chamber $Ch_2$ defined by 
\begin{center}
  $
\begin{cases}
a_1+a_2>1\\
a_1+a_3>1\\
a_1+a_4>1\\
a_2+a_3>1\\
a_2+a_4>1\\
a_3+a_4>1\\
a_i+a_j<1 \text{ for any other couple of indices}\\
\sum_{i\in S}a_i>1 \text{ only if } |S|\geq 3 \text{ and } |S\cap\{1,2,3,4\}|\geq 2.
\end{cases}
$  
\end{center}
This also is not empty since $$\mathcal{A}_2=(\frac{1}{2}+\varepsilon,\frac{1}{2}+\varepsilon,\frac{1}{2}+\varepsilon,\frac{1}{2}+\varepsilon,\varepsilon,\varepsilon,\varepsilon,\varepsilon)$$ belongs to $Ch_2$ for a sufficiently small $\varepsilon$.  We can see that $Ch_1$ and $Ch_2$ are not comparable with respect to the partial order on $\mathbf{K}$ by looking at their defining inequalities. Indeed, by definition, given $Ch_a$ and $Ch_b$ in $\mathbf{K}$ we say that $Ch_a\leq Ch_b$ if they are equal or, informally, all of their defining inequalities which have different direction are such that have $\sum_{i\in S}a_i<1$ for every $\mathcal{A}\in Ch_a$ and $\sum_{i\in S}b_i>1$ for every $\mathcal{B}\in Ch_b$.  But here we have that in $Ch_1$, $a_1+a_5>1$ while $a_1+a_5<1$ in $Ch_2$. Meanwhile, $a_2+a_4<1$ in $Ch_1,$ but $a_2+a_4>1$ in $Ch_2$ so the above condition is not satisfied. 

Moreover these two chambers are in different orbits of the action of $S_8$ on $\mathbf{K}$ since the inequalities with two variables and right direction defining $Ch_1$ all contains the variable $a_1$, while in $Ch_2$ there is no variable repeated in all the inequalities with two variables and right direction. So there is no permutation $\sigma$ in $S_8$ such that applied to $\sigma (Ch_1)=Ch_2$, and this implies that $[Ch_1]$ and $[Ch_2]$ are not comparable in $[\mathbf{K}]$.

At the level of moduli spaces of tropical curves, it means that we can find points in $M_{g,\mathcal{A}_1}^{trop}$ with combinatorial types which are not in $M_{g,\mathcal{A}_2}^{trop}$ and viceversa as shown in the Figure \ref{comparisonExample}.
\begin{center}
\begin{figure}[h]
\scalebox{0.6}{
\begin{tikzpicture}
\draw[black](0,0)--(1,1);
\draw[black](2,0)--(1,1);
\draw[black](1,1)--(1,2);
\fill[black] (1,2) circle (0.4ex);
\fill[black] (1,1) circle (0.4ex);
\fill[black] (2,0) circle (0.4ex);
\fill[black] (0,0) circle (0.4ex);
\draw[black](0,0)--(0,-0.5);
\draw[black](0,0)--(-0.5,-0.5);
\draw[black](0,0)--(-0.5,0);
\draw[black](2,0)--(2,-0.5);
\draw[black](2,0)--(2.5,-0.5);
\draw[black](2,0)--(2.5,0);
\draw[black](1,2)--(0.5,2.5);
\draw[black](1,2)--(1.5,2.5);

\node at (0.5,2.7){1};
\node at (1.5,2.7){2};
\node at (-0.8,0){3};
\node at (-0.7,-0.7){4};
\node at (0,-0.7){5};
\node at (2.8,0){6};
\node at (2.7,-0.7){7};
\node at (2,-0.7){8};
\node at (1,0.7){\textbf{g}};
\end{tikzpicture}

\qquad

\begin{tikzpicture}
\draw[black](0,0)--(2,0);
\fill[black] (0,0) circle (0.4ex);
\fill[black] (1,0) circle (0.4ex);
\fill[black] (2,0) circle (0.4ex);
\draw[black](0,0)--(-0.5,-0.5);
\draw[black](0,0)--(-0.5,0.5);
\draw[black](2,0)--(2.5,0.5);
\draw[black](2,0)--(2.5,-0.5);
\draw[black](1,0)--(0.3,0.3);
\draw[black](1,0)--(1.8,0.3);
\draw[black](1,0)--(0.75,0.5);
\draw[black](1,0)--(1.25,0.5);

\node at (-0.6,0.7){1};
\node at (-0.6,-0.7){2};
\node at (2.6,0.7){3};
\node at (2.6,-0.7){4};

\node at (0.3,0.5){5};
\node at (0.75,0.7){6};
\node at (1.25,0.7){7};
\node at (1.8,0.5){8};
\node at (1,-0.3){\textbf{g}};
\end{tikzpicture}

}
    \caption{The graph on the left is a combinatorial type of points in $M_{g,\mathcal{A}_1}^{trop}$ but not in $M_{g,\mathcal{A}_2}^{trop}$, while the graph on the right is a combinatorial type of points in $M_{g,\mathcal{A}_2}^{trop}$ but not in $M_{g,\mathcal{A}_1}^{trop}$.}\label{comparisonExample}
\end{figure}
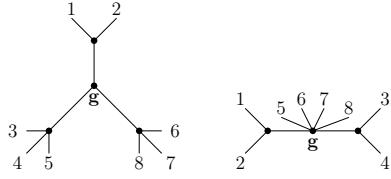
\end{center}
\end{myex}

\subsection{Wall crossing on graph complexes} Consider now the graph complexes $G^{(g,\mathcal{A})}$ introduced in Section \ref{complexes}. These complexes are defined upon the same stability conditions that we have been considering so far on graphs, so we can establish a theorem analogous to Theorem \ref{Main1} which holds for $G^{(g,\mathcal{A})}$.
\begin{myprop}\label{Compprop}
Let $g,n\geq 1$ be two integers, and $\mathcal{A}$ and $\mathcal{B}$ two weight data in $\mathcal{D}_{g,n}$.
\begin{itemize}
    \item[1)] If $\mathcal{A}\leq \mathcal{B}$, then $G^{(g,\mathcal{A})}\subset G^{(g,\mathcal{B})}$;
    \item[2)] If $\mathcal{A}$ and $\mathcal{B}$ are in the same chamber, then $G^{(g,\mathcal{A})}= G^{(g,\mathcal{B})}$;
    \item[3)] If $\mathcal{A}$ and $\mathcal{B}$ are obtained one from the other through a permutation of coordinates, then $G^{(g,\mathcal{A})}$ is isomorphic to $G^{(g,\mathcal{B})}$;
    \item[4)]If $\mathcal{A}$ and $\mathcal{B}$ are in chambers $Ch_1$ and $Ch_2$ such that $[Ch_1]=[Ch_2]\in[\mathbf{K}]$, then $G^{(g,\mathcal{A})}$ is isomorphic to $G^{(g,\mathcal{B})}$.
\end{itemize}
\end{myprop}
\begin{proof} A generator of  $G^{(g,\mathcal{A})}$ is a $(g,\mathcal{A})$-stable graph  with an edge ordering under the relation $[\G,\omega]=sgn(\sigma)[\G',\omega']$ if there is an isomorphism of $n$-marked graphs $\G\cong\G'$ under which the orders $\omega$ and $\omega'$ are related by the permutation $\sigma\in S_{|E(\G)|}$
If $\mathcal{A}\leq\mathcal{B}$ , every $(g,\mathcal{A})$-stable graph is also $(g,\mathcal{B})$-stable, so to show 1 we just send a generator into itself. For 2, if $\mathcal{A}$ and $\mathcal{B}$ are in the same chamber then $\mathcal{G}_{g,\mathcal{A}}=\mathcal{G}_{g,\mathcal{B}}$ so the result follows. Proof of 3 follows by the same reasoning by sending a generator of $G^{(g,\mathcal{A})}$ into the generator of $G^{(g,\mathcal{B})}$ obtained relabeling the legs according to the permutation which sends $\mathcal{A}$ to $\mathcal{B}$. This also shows 4 as by hypothesis we can assume $\mathcal{A}$ and $\mathcal{B}$ are obtained one from the other after a permutation of coordinates.
\end{proof}
\begin{myth}\label{Complexfiltration}
Let $g\geq 0$, $n\geq 1$ be two integers. Fix a weight datum $\mathcal{A}\in \mathcal{D}_{g,n}.$
There are filtrations of $G^{(g,\mathcal{A})}$ induced by the partial order on the set of chambers up to symmetry given by inclusions of complexes. Namely a ordered sequence $$[Ch_{\mathcal{A}_1}]\leq [Ch_{\mathcal{A}_2}]\leq ...\leq [Ch_{\mathcal{A}_p}] \leq ... \leq [Ch_{\mathcal{A}_{N-1}}]\leq [Ch_{\mathcal{A}}], $$ induces a filtration of chain complexes $$G^{(g,\mathcal{A}_1)}\hookrightarrow G^{(g,\mathcal{A}_2)}\hookrightarrow ...\hookrightarrow G^{(g,\mathcal{A}_p)} \hookrightarrow ... \hookrightarrow G^{(g,\mathcal{A}_{N-1})}\hookrightarrow G^{(g,\mathcal{A})}, $$
with
$G^{(g,\mathcal{A}_{p-1})}\hookrightarrow G^{(g,\mathcal{A}_p)}$ being an injective map of chain complexes for every $p=2,...,N$.
\end{myth}
\begin{proof}
 At each step of the filtration $[Ch_{\mathcal{A}_p}]\leq [Ch_{\mathcal{A}_{p+1}}]$ we can find two weight data $\mathcal{A}\leq \mathcal{B}$ and two chambers $Ch_1\leq Ch_2$ such that $\mathcal{A}\in Ch_1\in [Ch_{\mathcal{A}_p}]$ and $\mathcal{B}\in Ch_2\in [Ch_{\mathcal{A}_{p+1}}]$. Then by Proposition \ref{Compprop} $G^{(g,\mathcal{A}_p}) \cong G^{(g,\mathcal{A})}\subset G^{(g,\mathcal{B})} \cong G^{(g,\mathcal{A}_{p+1})}$ analogously to what we did to prove Theorem \ref{Main1}.  
\end{proof}
\subsection{Examples of filtrations}
In this section we construct some sequences of weight data, chambers and chambers up to symmetry which define interesting filtrations of moduli spaces. 
\subsection*{The five term filtration} 
 Recall that for any $g,n\geq1$ there are two special chambers up to symmetry: the maximal chamber $[Ch_{1^{(n)}}]$, which is greater  than any other chamber up to symmetry, and the minimal chamber $[Ch_{\varepsilon^{(n)}}]:=[Ch_{\mathcal{E}}]$, for $\varepsilon<\frac{1}{n}$, which is smaller than any other chamber up to symmetry. Moreover, both these orbits are made by a single chamber. 

Let $\mathcal{F}=(\frac{1}{n}+\varepsilon,\frac{1}{n},...,\frac{1}{n})$. This datum belongs to the chamber $Ch_{\mathcal{F}}$ described by the inequalities $\sum_{i\in S}a_i<1$ for every $S\neq \{1,...,n\}$ and $\sum_{i=1}^na_i>1$ . It is clearly invariant by the action of $S_n$ and so its orbit $[Ch_{\mathcal{F}}]$ is made only by itself. In particular, $\sigma(\mathcal{F})\in Ch_{\mathcal{F}}$ for every $\sigma\in S_n$. Moreover, by the inequalities defining $Ch_{\mathcal{F}}$ it is clear that $Ch_{\varepsilon^{(n)}}\leq Ch_{\mathcal{F}}\leq Ch$ for any chamber $Ch\neq Ch_{\varepsilon^{(n)}}$.

Suppose we have two chambers up to symmetry $[Ch_1]$ and $[Ch_2]$ such that $[Ch_1]\leq [Ch_2]$: then it is always possible to construct a five term sequence
$$[Ch_{\mathcal{E}}]\leq [Ch_{\mathcal{F}}]\leq [Ch_1]\leq [Ch_2] \leq [Ch_{1^{(n)}}]$$ inducing a filtration of moduli spaces
$$M_{g,\mathcal{E}}^{trop}\hookrightarrow M_{g,\mathcal{F}}^{trop}\hookrightarrow M_{g,\mathcal{A}}^{trop}\hookrightarrow M_{g,\mathcal{B}}^{trop}\hookrightarrow M_{g,n}^{trop},$$ where $\mathcal{A}\in Ch_1$ and $\mathcal{B}\in Ch_2$.

\begin{myrem}
In general, given a sequence $[Ch_{\mathcal{A}_1}]\leq...\leq [Ch_{\mathcal{A}_p}]$ which does not contain already $[Ch_{\mathcal{E}}]$, $[Ch_{\mathcal{F}}]$ and $[Ch_{1^(n)}$ , it is always possible to extend it by two terms $[Ch_{\mathcal{E}}]\leq [Ch_{\mathcal{F}}]$ on the left and by $[Ch_{1^{(n)}}]$ on the right (if they are not already in the sequence).
\end{myrem}
\subsection*{The Heavy/Light filtration}
Let $g\geq 1$, $n\geq 2$, and $\varepsilon<\frac{1}{n}$. We denote by $(1^{(m)},\varepsilon^{(n-m)})$ the weight datum made by a sequence of $m$ components of weight 1, called the heavy components, and $n-m$ components of weight $\varepsilon^{(n-m)}$, called light components, for $m$ which goes from $0$ to $n$. Such data are called heavy/light. Notice that if $m=0$ we get $\mathcal{E}$, while if $m=n$ we get $1^{(n)}$. 

By construction, we have $$(1^{(m)},\varepsilon^{(n-m)})\leq(1^{(m+1)},\varepsilon^{(n-m-1)})$$ for any $m=0,...,n-1$, and passing to their chambers it is easy to see that we have
$$Ch_{\mathcal{E}}\leq Ch_{(1^{(1)},\varepsilon^{(n-1)})}\leq...\leq Ch_{(1^{(m)},\varepsilon^{(n-m)})}\leq ... \leq Ch_{(1^{(n-1)},\varepsilon^{(1)})}\leq 1^{(n)}.$$
It is easy to check that all of these inequalities are strict except the last one, since $Ch_{(1^{(n-1)},\varepsilon^{(1)})}=Ch_{ 1^{(n)}}$. 

So there is an induced filtration of spaces
$$M_{g,\mathcal{E}}^{trop}\hookrightarrow M_{g,(1^{(1)},\varepsilon^{(n-1)})}^{trop}\hookrightarrow...\hookrightarrow M_{g,(1^{(m)},\varepsilon^{(n-m)})}^{trop}\hookrightarrow...\hookrightarrow M_{g,n}^{trop}.$$

Of course the same holds also if we consider the chambers up to symmetry. 
In that case, given a weight datum of the heavy/light form $(1^{(m)},\varepsilon^{(n-m)})$, each chamber in $[Ch_{(1^{(m)},\varepsilon^{(n-m)})}]$ contains a heavy/light datum obtained by a permutation of the original one.
\begin{myrem}
The heavy/light filtration also works for $g=0$, but in this case we must start by the case $m=2$. These particular input data were used in \cite{CHMR2014moduli} and \cite{kannan2021chow} to study the topology and the Chow Ring of the moduli spaces of rational tropical curves.
\end{myrem}
\subsection*{The floor filtration}\label{floor}
Fix $g\geq 0$ and $l\in\{2,...,n\}$. Suppose we have a chamber defined by the following set of inequalities: $\sum_{i\in S}a_i>1$ if and only if $|S|\geq l$. The weight datum $$\mathcal{H}_l:=(\frac{1}{l}+\varepsilon_l,...,\frac{1}{l}+\varepsilon_l)$$ belongs to that chamber for $\varepsilon_l$ sufficiently small, and its chamber $Ch_{\mathcal{H}_l}$ is invariant with respect to the action of $S_n$. We call $[Ch_{\mathcal{H}_l}]$ the $l$-th floor of the chamber decomposition.

Notice that if $l=2$, $Ch_{\mathcal{H}_2}=Ch_{1^{(n)}}$, while if $l=n$,  $Ch_{\mathcal{H}_n}=Ch_{\mathcal{F}}$. Moreover by construction $$Ch_{\mathcal{H}_l}\leq Ch_{\mathcal{H}_{l-1}}$$
(and clearly the same holds if we pick them up to symmetry) so there is an induced filtration
$$M_{g,\mathcal{F}}^{trop}\hookrightarrow...\hookrightarrow M_{g,\mathcal{H}_l}^{trop}\hookrightarrow M_{g,\mathcal{H}_{l-1}}^{trop}...\hookrightarrow  M_{g,n}^{trop},$$
eventually extendable on the left with $M_{g,\mathcal{E}}^{trop}$.
Notice that being stable for a curve in the $l-$th floor means that each valence one vertex needs at least $l$ markings to be stable.
\section{Spectral sequences of graph complexes}\label{section6}

We start this section putting in relation the reduced rational homology of $\Delta_{g,\mathcal{A}}$ with the top weight cohomology of $\mathcal{M}_{g,\mathcal{A}}$. We start giving a brief account on  the theory of symmetric $\Delta$-complexes introduced in  \cite{CGP2} and \cite{CGP1}.

\subsection*{Symmetric $\Delta$-complexes} Let $I$ denote the category whose
objects are the sets $[p]:=\{0,...,p\}$ for non-negative integers $p$, together with $[-1]:=\emptyset$,
and whose morphisms are injections of sets. A symmetric $\Delta$-complex is a functor from $I^{op}$ to
$Sets$.

For any weight datum $\mathcal{A}$, we can consider our $\Delta_{g,\mathcal{A}}$'s as symmetric $\Delta$-complexes as follows. Let $X=\Delta_{g,\mathcal{A}}:I^{op}\rightarrow Sets$ be a functor, with $$X_p=\{\text{equivalence classes of pairs } (G,\tau)\},$$ with $G\in Ob(\mathcal{G}_{g,\mathcal{A}})$ and $\tau:E(G)\rightarrow [p]$ a bijection, where we consider $\tau=\tau'$ if they are in the same orbit under the action of $Aut(G)$. For every $i:[p']\rightarrow [p]$ define $X(i):X_p\rightarrow X_{p'}$ as follows: given an element of $X_p$ represented by $(G,\tau:E(G)\rightarrow [p])$ we contract the edges of $G$ whose labels are not in $i([p'])\subset [p]$, and then we relabel the remaining edges with labels $[p']$ as prescribed by $i$. The result is a $[p']$ edge labeled graph $G'$, and we set it to be $X(i)(G)$. 

\noindent To a symmetric $\Delta$-complex $X$, we can associate its group of cellular $p$-chains $$C_p(X)=(\mathbb{Q}^{sign}\otimes \mathbb{Q}X_p)_{S_{p+1}}$$ where $\mathbb{Q}X_p$ is the vector space with basis $X_p$ on which $S_{p+1}$ acts by permuting the basis vectors, and $\mathbb{Q}^{sgn}$ denotes the action of $S_{p+1}$ on $\mathbb{Q}$ via the sign. By Proposition 3.8 of \cite{CGP1} we know that the homology of $C_*(\Delta_{g,\mathcal{A}})$ is identified with $\widetilde{H}_*(\Delta_{g,\mathcal{A}};\mathbb{Q})$.

Moreover, whenever $X\subset Y$ is a subcomplex in the sense of Definition 3.5 of \cite{CGP2}, for every $p\geq -1$ we can define $C_p(Y,X)$ by the exact sequence $$0 \rightarrow C_p(X)\rightarrow C_p(Y)\rightarrow C_p(Y,X)\rightarrow 0.$$ 

For every weight data $\mathcal{A}\in Ch_1$ and $\mathcal{B}\in Ch_2$ such that $[Ch_1]\leq [Ch_2]$, since there is an inclusion as symmetric $\Delta$-subcomplex $\Delta_{g,\mathcal{A}}\subset \Delta_{g,\mathcal{B}}$, we can consider the relative chain complex $C_*(\Delta_{g,\mathcal{B}},\Delta_{g,\mathcal{A}})$ whose homology is identified with the relative rational homology by Proposition 3.6., \cite{CGP2}:
$$H_i(C_*(\Delta_{g,\mathcal{B}},\Delta_{g,\mathcal{A}}))\cong H_i(\Delta_{g,\mathcal{B}},\Delta_{g,\mathcal{A}};\mathbb{Q}).$$

First, we can prove the following Theorem. 
\begin{myth}\label{HomTh}
Let $g\geq 1$ and $\mathcal{A}\in \mathcal{D}_{g,n}$ a weight datum. There is a natural injection of chain complexes $$G^{(g,\mathcal{A})}\rightarrow C_*(\Delta_{g,\mathcal{A}},\mathbb{Q}) $$ decreasing degrees by $2g-1$, inducing isomorphisms on homology $$\widetilde{H}_{k+2g-1}(\Delta_{g,\mathcal{A}};\mathbb{Q})\rightarrow H_k(G^{(g,\mathcal{A})})$$ for all $k$'s.
\end{myth}
\begin{proof}
We prove this theorem following the lines of Theorem 1.4 of \cite{CGP2}.
Consider the cellular chain complex $C_*(\Delta_{g,\mathcal{A}}, \mathbb{Q})$. It is generated in degree $p$ by
$[G,\omega]$ where $G\in\mathcal{G}_{g,\mathcal{A}}$ is a graph and $\omega : E(G)\rightarrow [p] = {0, 1,..., p}$ is a bijection, with the relations $[G,\omega]= sgn(\sigma)[G',\omega']$ if there is an isomorphism $G\rightarrow G'$ of graphs 
inducing the permutation $\sigma$ of the set $[p]$.
We claim that the complex $C_*(\Delta_{g,\mathcal{A}}, \mathbb{Q})$ splits into the direct sum of two subcomplexes $A^{(g,\mathcal{A})}\bigoplus B^{(g,\mathcal{A})}$, $A^{(g,\mathcal{A})}$ being spanned by the generators  where $G$ has no loops and whose vertices have weight zero, and $B^{(g,\mathcal{A})}$ is spanned by the generators such that $G$ has at least one loop or one nonzero vertex weight. These are in fact subcomplexes: for $B^{(g,\mathcal{A})}$ it is clear, while   for $A^{(g,\mathcal{A})}$, we need just to observe that if
$G$ has no loops and all vertices have weight zero, and $[G,\omega]\neq 0$, then $G$ has no parallel
edges. Therefore, every
contraction $G/e$ also has no loops and also has all vertices of weight zero.

Now we note that $A^{(g,\mathcal{A})}$ is isomorphic to  $G^{(g,\mathcal{A})}$ up to
shifting degrees by $2g - 1$, and $B^{(g,\mathcal{A})}$ is the cellular chain complex associated to the subcomplex of tropical curves with underlying graphs having at least a loop and/or a vertex with positive weight $\Delta^{lw}_{g,\mathcal{A}}$, which is contractible whenever it is nonempty by Theorem 3.2 of \cite{KLSY}. Therefore $B^{(g,\mathcal{A})}$
is an
acyclic complex by Proposition 3.8 of \cite{CGP1}, so the result follows.
\end{proof}
So the graph complexes $G^{(g,\mathcal{A})}$ compute the reduced rational homology of the $\Delta_{g,\mathcal{A}}$. This Theorem generalizes Theorem 1.4 of \cite{CGP2}, and gives the following Corollary.
\begin{mycor}\label{TopW}
There is a natural isomorphism $$Gr^W_{6g-6+2n}H^{4g-6+2n-k}(\mathcal{M}_{g,\mathcal{A}};\mathbb{Q})\rightarrow H_{k}(G^{(g,\mathcal{A})}).$$
\end{mycor}
\begin{proof}
Using the language of \cite{CGP1} and \cite{CGP2}, the proof follows by Theorem \ref{HomTh} and Theorem 5.8 of \cite{CGP1} just observing that the dual boundary complex $\Delta(\mathcal{M}_{g,\mathcal{A}}\subset\overline{\mathcal{M}}_{g,\mathcal{A}})$ is $\Delta_{g,\mathcal{A}}.$ Indeed, by Theorem 1.1 of \cite{ulirsch14tropical}, we have $D:=\overline{\mathcal{M}}_{g,\mathcal{A}}\setminus \mathcal{M}_{g,\mathcal{A}}$ to be a divisor with normal crossing (stack theoretically), so by Corollary 5.6 of \cite{CGP1}, $\Delta(D)$ is the symmetric $\Delta$-complex associated to the smooth generalized cone complex $\mathfrak{S}(\overline{\mathcal{M}}_{g,\mathcal{A}})$. The rest of the proof is analogous to the one of Corollary 5.6 of \cite{CGP1} and Theorem 6.1 of \cite{CGP2}.
\end{proof}

\subsection*{Filtered Chain Complexes}
We want to put on the graph complexes the structure of filtered chain complexes and extract from them a spectral sequence in order to show our last theorem. We recall briefly some definitions and facts about these objects.
\begin{mydef}\label{filcom}
A filtered module is an $R$-module $A$ with an increasing sequence of sub-modules $F_pA\subset F_{p+1}A$ indexed by $p\in \mathbb{Z}$ such that $\cup_{p\in \mathbb{Z}} F_pA=A$ and $\cap_{p\in \mathbb{Z}}F_pA=\{0\}$. We call $\{F_pA\}_{p\in \mathbb{Z}}$ a filtration of $A$.
\end{mydef}
We say that the filtration $\{F_pA\}_{p\in \mathbb{Z}}$ is bounded if $F_pA=\{0\}$ for sufficiently small $p$ and $F_pA=A$ for sufficiently large $p$.
\begin{mydef}
Let $A$ be a filtered module. The associated graded module is defined, in degree $p$, as $G_pA=F_pA/F_{p-1}A$.
\end{mydef}
Notice that there is a short exact sequence $$0\rightarrow F_{p-1}A\rightarrow F_pA\rightarrow G_pA\rightarrow 0.$$
\begin{mydef}
A filtered chain complex is a chain complex $(C_*,\partial)$ together with a filtration $\{F_pC_i\}_{p\in \mathbb{Z}}$ on each $C_i$ such that the differential preserves the filtration, i.e. $\partial(F_pC_i)\subset F_pC_{i-1}$.
\end{mydef}
We have a well defined induced differential $\partial:G_pC_i\rightarrow G_pC_{i-1}$, and so we can define an associated graded chain complex $G_pC_*.$ Moreover, there is an induced filtration on the homology of $C_*$ given by
$$F_pH_i(C_*)=\{\alpha\in H_i(C_*)|\alpha =[x],\exists x\in F_pC_i\}.$$
Again, this has associated graded pieces $G_pH_i(C_*)$ defined as before.
\begin{mydef}
A spectral sequence consists of:
\begin{itemize}
    \item An $R$-module $E^r_{p,q}$ defined for each $p,q \in \mathbb{Z}$ and each integer $r\geq r_0$,
    where $r_0$ is some nonnegative integer;
    \item Differentials $\partial_r : E^r_{p,q} \rightarrow E^r_{p-r,q+r-1}$ such that $\partial_r^2 = 0$ and $E^{r+1}$ is the homology of $(E^r,\partial_r)$, i.e.
$$E^{r+1}_{p,q} = \frac{Ker(\partial_r : E^r_{p,q} \rightarrow
E^1_{p-r,q+r-1})}{Im(\partial_r : E^1_{p+r,q-r+1} \rightarrow E^1_{p,q})} 
 $$
\end{itemize}
\end{mydef}

A spectral sequence converges if for every $p, q$, if $r$ is sufficiently large then $\partial_r$  vanishes on $E^p_{r,q}$ and $E^r_{p+r,q-r+1}$. In this case, for each $p,q$, the module
$E^r_{p,q}$ is independent of $r$ for $r$ sufficiently large, and we denote this by $E^{\infty}_{p,q}$.
For a given $r$, the collection of $R$-modules $\{ E^r_{p,q} \}$, together with the differentials $\partial_r$ between them, is called the $r$-th page of the spectral sequence. Each page is the homology of the previous page. 

Given a filtered complex, we have an associated spectral sequence obtained from the short exact sequences extracted from the filtrations. Namely, let $E^0_{p,q}:=G_p C_{p+q}$. Then there is a well defined $\partial:E^0_{p,q}\rightarrow E^0_{p,q-1}.$ We denote $E^1_{p,q} =H_{p+q}(G_pC_*).$
and define $\partial_1:E^1_{p,q}\rightarrow E^1_{p-1,q}$ as follows. A homology class $\alpha \in E^1_{p,q}$ can be represented by a chain $x\in F_pC_{p+q}$ such that $\partial x \in F_{p-1}C_{p+q-1}$. We set $\partial_1(\alpha) = [\partial x]$. It follows easily from $\partial^2_1 = 0$ that $\partial_1$ is well defined and $\partial_1^2 = 0$. We now consider the homology
$$E^2_{p,q} =\frac{Ker(\partial_1 : E^1_{p,q} \rightarrow E^1_{p-1,q})}{Im(\partial_1 : E^1_{p+1,q} \rightarrow E^1_{p,q})}.$$
We can iterate this process for every nonnegative integer $r$, so we can define an $r$-th order approximation to $G_pH_{p+q}(C_*)$ by $$E^r_{p,q} = \frac{\{x \in F_pC_{p+q} | \partial x \in F_{p-r}C_{p+q-1} \}}{F_{p-1}C_{p+q} + \partial(F_{p+r-1}C_{p+q+1})}.$$
 
Here the notation indicates the quotient of the numerator by its intersection with the denominator.

\noindent Fix $g,n\geq 1$ and a 
sequence $[Ch_{\mathcal{A}_1}]\leq...\leq [Ch_{\mathcal{A}_N}]$. 
It gives a filtration of chain complexes
$$G^{(g,\mathcal{A}_1)}\hookrightarrow...\hookrightarrow G^{(g,\mathcal{A}_p)} \hookrightarrow ... \hookrightarrow G^{(g,\mathcal{A}_N)};$$

\noindent If we set by convention $G^{(g,\mathcal{A}_p)}=\{0\}$ for every $p\leq 0$ and we let it stabilize at the last term for every $p\geq N$, we can extend the above filtration for every $p\in \mathbb{Z}$, with the trivial differential outside the bounds 0 and $N$. Then the induced filtration on each $G^{(g,\mathcal{A}_N)}_i$ makes $G^{(g,\mathcal{A}_N)}$ a filtered chain complex. 

\noindent This structure on $G^{(g,\mathcal{A}_N)}$ induces a spectral sequence as already seen, and we have $$E^0_{p,q}=G_pG^{(g,\mathcal{A}_N)}_{p+q}=F_pG^{(g,\mathcal{A}_N)}_{p+q}/F_{p-1}G^{(g,\mathcal{A}_N)}_{p+q}=G_{p+q}^{(g,\mathcal{A}_p)}/G_{p+q}^{(g,\mathcal{A}_{p-1})}:=G_{p+q}^{(g,\mathcal{A}_p,\mathcal{A}_{p-1})},$$
with $G_{p+q}^{(g,\mathcal{A}_p,\mathcal{A}_{p-1})}$ being the complex of $\mathcal{A}_p$ but not $\mathcal{A}_{p-1}$ stable graphs, and
$$E^1_{p,q}=H_{p+q}(G_pG^{(g,\mathcal{A}_N)})=H_{p+q}(G_{p+q}^{(g,\mathcal{A}_p)}/G_{p+q}^{(g,\mathcal{A}_{p-1})})=H_{p+q}(G_{p+q}^{(g,\mathcal{A}_p,\mathcal{A}_{p-1})}).$$
So we are ready to show the last main Theorem, restated here.

\begin{theoremb*}
Fix $g\geq 1$, $n\geq 2$. Assume we have a sequence of chambers up to symmetry $[Ch_{\mathcal{A}_1}]\leq...\leq [Ch_{\mathcal{A}_p}]\leq...\leq [Ch_{\mathcal{A}_N}]$, and let $G^{(g,\mathcal{A}_1)}\hookrightarrow...\hookrightarrow G^{(g,\mathcal{A}_p)}\hookrightarrow...\hookrightarrow G^{(g,\mathcal{A}_N)}$ be the induced filtration on graph complexes. Then 
$$Gr^W_{6g-6+2n}H^{4g-6+2n-k}(\mathcal{M}_{g,\mathcal{A}_N};\mathbb{Q})\cong \bigoplus_{p=1}^NE^{\infty}_{p,k-p},$$
where the terms $E^{\infty}_{p,k-p}$ are the ones to which the spectral sequence induced by the above filtration converges.
\end{theoremb*}
\begin{proof}
We consider the homology ring $H_*(G^{(g,\mathcal{A}_N)})=\bigoplus_{k\in\mathbb{Z}}H_k(G^{(g,\mathcal{A}_N)})$. By construction of the graph complexes $H_*(G^{(g,\mathcal{A}_N)})$ is supported only in degrees $1-2g\leq k\leq g+n-3$ so the sum is finite, and since any of the $H_k(G^{(g,\mathcal{A}_N)})$'s is finite dimensional this is true also for $H_*(G^{(g,\mathcal{A}_N)})$. Moreover, each $H_k(G^{(g,\mathcal{A}_N)})$ comes with a filtration induced by the one on $G^{(g,\mathcal{A}_N)}$: $$F_pH_k(G^{(g,\mathcal{A}_N)})=\{\alpha\in H_k(G^{(g,\mathcal{A}_N)})|\alpha=[x],x\in G_k^{(g,\mathcal{A}_p)}\},$$
so also $H_*(G^{(g,\mathcal{A}_N)})$ is filtered by $$F_pH_*(G^{(g,\mathcal{A}_N)})=\bigoplus_{k=1-2g}^{g+n-3}F_pH_k(G^{(g,\mathcal{A}_N)}).$$
This gives to $H_*(G^{(g,\mathcal{A}_N)})$ the structure of a finite dimensional filtered graded vector space. 

By \cite{mccleary2001user}, Section 1 such an object can be decomposed, in each degree, as the direct sum $$H_k(G^{(g,\mathcal{A}_N)})=\bigoplus_{p+q=k}\frac{F_pH_{p+q}(G^{(g,\mathcal{A}_N)})}{F_{p-1}H_{p+q}(G^{(g,\mathcal{A}_N)})}=\bigoplus_{p+q=k}G_pH_{p+q}(G^{(g,\mathcal{A}_N)}).$$
Moreover, since the filtration is bounded, we can rewrite this sum taking only the significant indices:
$$H_k(G^{(g,\mathcal{A}_N)})=\bigoplus_{p=1}^N\frac{F_pH_{p+q}(G^{(g,\mathcal{A}_N)})}{F_{p-1}H_{p+q}(G^{(g,\mathcal{A}_N)})}=\bigoplus_{p=1}^NG_pH_{p+q}(G^{(g,\mathcal{A}_N)})).$$
where $q=k-p$.

Consider now the associated spectral sequence: by construction we have $$G_pH_{p+q}(G^{(g,\mathcal{A}_N)})=E^{\infty}_{p,q}$$ since the spectral sequence converges to it, so
$$H_k(G^{(g,\mathcal{A}_N)})=\bigoplus_{p=1}^NE^{\infty}_{p,k-p}.$$
To conclude it is enough to apply Corollary \ref{TopW}, so that $$Gr^W_{6g-6+2n}H^{4g-6+2n-k}(\mathcal{M}_{g,\mathcal{A}_N};\mathbb{Q})\cong H_k(G^{(g,\mathcal{A}_N)})= \bigoplus_{p=1}^NE^{\infty}_{p,k-p}.$$
\end{proof}
We add some remarks about the proof:
\begin{myrem}
The proof also shows that the same decomposition holds for $\widetilde{H}_{k-1}(\Delta_{g,\mathcal{A}_N};\mathbb{Q})$, just by applying Theorem \ref{HomTh} instead of the Corollary.
\end{myrem}
\begin{myrem}
Since the filtration is bounded, by Lemma 3.1.d of \cite{hutchings2011introduction} the approximations of the spectral sequence stabilize after a certain $r$, i.e $E^{\infty}_{p,q}=E^r_{p,q}$ for $r$ sufficiently large. By the description $$E^r_{p,q} = \frac{\{x \in F_pC_{p+q} | \partial x \in F_{p-r}C_{p+q-1} \}}{F_{p-1}C_{p+q} + \partial(F_{p+r-1}C_{p+q+1})}$$
we can see that for our filtrations these terms stabilize when $r\geq max\{p,N-p+1\}$ since for these $r$'s $F_{p-r}C_{p+q-1}=\{0\}$ and $F_{p+r-1}C_{p+q+1}=G^{(g,\mathcal{A}_N)}_{p+q}$. So $E^{\infty}_{p,q}=E^{max\{p,N-p+1\}}_{p,q}$.
\end{myrem}
\begin{myrem}
 The top weight cohomology $Gr^W_{6g-6+2n}H^{4g-6+2n-k}(\mathcal{M}_{g,\mathcal{A}_N})$ does not depend on the chosen filtration: a priori different filtrations give different spectral sequences and different convergence terms $E^{\infty}_{p,k-p}$, Nevertheless, the two direct sums will be isomorphic.
 
\noindent In general, some of the terms $E^{\infty}_{p,k-p}$ in the direct sum may be zero.
\end{myrem}

\begin{myex}\label{FtComp}
We put ourselves in the case $g=1, n=3$. We fix the sequence of weight data of Example \ref{final}: $$ (\frac{1}{3},\frac{1}{3},\frac{1}{3}-\varepsilon) \leq (\frac{4}{9}-\varepsilon,\frac{4}{9}-\varepsilon,\frac{4}{9}-\varepsilon) \leq (\frac{14}{27}-\varepsilon, \frac{12}{27}, \frac{14}{27}) \leq (1-\varepsilon, \frac{12}{27},\frac{14}{27}) \leq 1^{(n)},$$
for a sufficiently small $\varepsilon$. 
We set also $\mathcal{A}_1=(\frac{1}{3},\frac{1}{3},\frac{1}{3}-\varepsilon),$  $\mathcal{A}_2=(\frac{4}{9}-\varepsilon,\frac{4}{9}-\varepsilon,\frac{4}{9}-\varepsilon)$,  $\mathcal{A}_3=(\frac{14}{27}-\varepsilon, \frac{12}{27}, \frac{14}{27})$,  $\mathcal{A}_4=(1-\varepsilon, \frac{12}{27},\frac{14}{27})$ and $\mathcal{A}_5=1^{(n)}$  for simplicity.
The corresponding sequence of chambers up to symmetry is then
$$[Ch_{(\frac{1}{3},\frac{1}{3},\frac{1}{3}-\varepsilon) }]\leq[Ch_{(\frac{4}{9}-\varepsilon,\frac{4}{9}-\varepsilon,\frac{4}{9}-\varepsilon)}]\leq [Ch_{(\frac{14}{27}-\varepsilon, \frac{12}{27}, \frac{14}{27})}]\leq [Ch_{(1-\varepsilon, \frac{12}{27},\frac{14}{27})}]\leq[Ch_{1^{(n)}}].$$

By Theorem \ref{LastTheorem}, we can compute the homology of $\Delta_{1,1^{(n)}}=\Delta_{1,3}$, which corresponds to the top weight cohomology of $\mathcal{M}_{1,3}$, computing first the terms of the spectral sequence of $G^{(1,3)}$ and then using the shifting degree isomorphism. Since we know that $G^{(1,3)}$ has homology only in degrees -1, 0, and 1, we have to compute only three direct sums, each with five terms. All the terms $E^{r}_{p,k-p}$ stabilize for $r\geq 5$.
\begin{itemize}
    \item Case $k=-1$. Here we have to compute $\bigoplus_{p=1}^5E^{5}_{p,-1-p}$. For each $p$ from $1$ to $5$, the term $E^{5}_{p,-1-p}$ is equal to the quotient $G_{-1}^{(1,\mathcal{A}_p)}/(G_{-1}^{(1,\mathcal{A}_{p-1})}+\partial G_{0}^{(1,\mathcal{A}_{p-1})})$. Now for every $p=1,...5$, $G_{-1}^{(1,\mathcal{A}_p)}$ is generated by the loop graph $L$ with a single vertex and its only orientation depicted in Figure \ref{fig:zero-1eee}, so the quotient is zero whenever $p\geq2$.
    \begin{center}
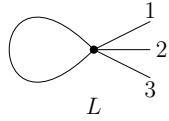
\begin{figure}[h]
\scalebox{0.75}{
\begin{tikzpicture}
\fill[black] (0,0) circle (0.4ex);
\draw[black] (0,0)--(1,0.5);
\draw[black] (0,0)--(1,-0.5);
\draw[black] (0,0)--(1,0);

\draw[black](0,0).. controls(-2,-2) and (-2,2) .. (0,0);

\node at (1,0.7){1};
\node at (1,-0.7){3};
\node at (1.2,0.){2};
\node at (0,-1){$L$};
\end{tikzpicture}}
\caption{This graph has only the loop-flipping automorphism, and only a possible order on its set of edges. The resulting generator $[L,\omega]$ has degree -1.}
\label{fig:zero-1eee}
\end{figure}
\end{center}
    When $p=1$, $\partial G_{0}^{(1,\mathcal{A}_{p-5+1})}=G_{-1}^{(1,\mathcal{A}_5)}=G_{-1}^{(1,\mathcal{A}_1)}$ by what we saw, so the quotient is again zero. Hence $$Gr^W_{6}H^{5}(\mathcal{M}_{1,3};\mathbb{Q})\cong\bigoplus_{p=1}^5E^{5}_{p,-1-p}=\{0\}.$$
    \item Case $k=0$: In this case, the terms are $E^{5}_{p,-p}$, for $p=1,...,5$. On the numerator of the quotient which defines $E^{5}_{p,-1-p}$ we have $Ker(\partial:G_0^{(1,\mathcal{A}_p)}\rightarrow G_{-1}^{(1,\mathcal{A}_p)})$ for each $p$ by construction. For $p=1$ this is zero since $G_0^{(1,\mathcal{A}_1)}$ is already zero. When $p=2$, $G_0^{(1,\mathcal{A}_2)}$ is generated by $[G,\omega_i]$, where $G$ is the graph of Figure \ref{fig:one-1ee2}.
    \begin{figure}[h]
\centering
\scalebox{0.75}{
\begin{tikzpicture}

\fill[black] (0,0-3) circle (0.4ex);
\fill[black] (-2,0-3) circle (0.4ex);
\draw[black] (0,0-3)--(1,0.5-3);
\draw[black] (0,0-3)--(1,-0.5-3);
\draw[black] (0,0-3)--(1,0-3);
\draw[black](-2,0-3).. controls(-4,-2-3) and (-4,2-3) .. (-2,0-3);
\draw[black](0,0-3)--(-2,0-3);

\node at (1,0.7-3){1};
\node at (1,-0.7-3){3};
\node at (1.2,0-3){2};
\node at (-1,0.3-3){$G$};

\end{tikzpicture}
}
\caption{The graph $G$ has no automorphisms, and since it has only two edges there are only two possible orderings. It gives two generators $[G,\omega_i]$ in degree 0, for $i=1,2$, related by their sign.}
\label{fig:one-1ee2}
\end{figure}
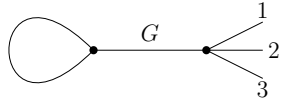

\noindent All the $\omega_i$ are related by permutations of the indices, so there is only a generator class.
The map $\partial$ sends the generators $[G,\omega_i]$ into $[L,\omega]$, so it is an isomorphism and $Ker(\partial:G_0^{(1,\mathcal{A}_2)}\rightarrow G_{-1}^{(1,\mathcal{A}_2)})$ is trivial. When $p=3,4,5$, the kernel $Ker(\partial:G_0^{(1,\mathcal{A}_p)}\rightarrow G_{-1}^{(1,\mathcal{A}_p)})$ coincides with $Im(\partial:G_0^{(1,\mathcal{A}_p)}\rightarrow G_{-1}^{(1,\mathcal{A}_p)})$, due to the generators coming from the graphs $G_1,$ $G_2$ and $G_3$ of figure \ref{fig:two-1ee}.
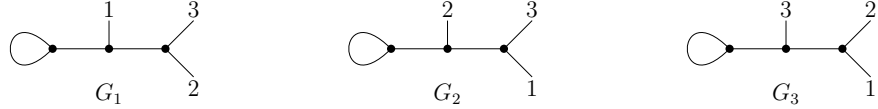
\begin{figure}[h]
\centering
\scalebox{0.75}{
\begin{tikzpicture}
\draw[black](2-6,0)--(1-6,0);
\draw[black](0-6,0).. controls(-1-6,1) and (-1-6,-1) .. (0-6,0);
\draw[black](0-6,0)--(1-6,0);

\draw[black](1-6,0)--(1-6,0.5);
\draw[black](2-6,0)--(2.5-6,-0.5);
\draw[black](2-6,0)--(2.5-6,0.5);

\fill[black] (0-6,0) circle (0.4ex);
\fill[black] (1-6,0) circle (0.4ex);
\fill[black] (2-6,0) circle (0.4ex);

\node at (1-6,0.7){$1$};
\node at (2.5-6,-0.7){$2$};
\node at (2.5-6,0.7){$3$};
\node at (1-6,-0.8){$G_1$};

\draw[black](2,0)--(1,0);
\draw[black](0,0).. controls(-1,1) and (-1,-1) .. (0,0);
\draw[black](0,0)--(1,0);

\draw[black](1,0)--(1,0.5);
\draw[black](2,0)--(2.5,-0.5);
\draw[black](2,0)--(2.5,0.5);

\fill[black] (0,0) circle (0.4ex);
\fill[black] (1,0) circle (0.4ex);
\fill[black] (2,0) circle (0.4ex);

\node at (1,0.7){$2$};
\node at (2.5,-0.7){$1$};
\node at (2.5,0.7){$3$};
\node at (1,-0.8){$G_2$};

\draw[black](2+6,0)--(1+6,0);
\draw[black](0+6,0).. controls(-1+6,1) and (-1+6,-1) .. (0+6,0);
\draw[black](0+6,0)--(1+6,0);

\draw[black](1+6,0)--(1+6,0.5);
\draw[black](2+6,0)--(2.5+6,-0.5);
\draw[black](2+6,0)--(2.5+6,0.5);

\fill[black] (0+6,0) circle (0.4ex);
\fill[black] (1+6,0) circle (0.4ex);
\fill[black] (2+6,0) circle (0.4ex);

\node at (1+6,0.7){$3$};
\node at (2.5+6,-0.7){$1$};
\node at (2.5+6,0.7){$2$};
\node at (1+6,-0.8){$G_3$};
\end{tikzpicture}}
\caption{The graphs with this combinatorial type have only an automorphism which flips the loop. This does not produce any odd permutation on the set of edges, so these classes are non-zero. Again there are six possible orderings on the edges, so we obtain generators  $[G_j,\omega_{j,i}]$ in degree 1, for $i=1,...,6$ and $j=1,2,3$.}
\label{fig:two-1ee}
\end{figure}
But $Im(\partial:G_0^{(1,\mathcal{A}_p))}\rightarrow G_{-1}^{(1,\mathcal{A}_p)})$ is exactly what we have in the denominator of the quotient defining $E^5_{p,-p}$, so it is zero. Hence all the terms of the direct sum are zero, so
$$Gr^W_{6}H^{4}(\mathcal{M}_{1,3};\mathbb{Q})=\{0\}.$$
 \item Case $k=1$: Here we have to compute $\bigoplus_{p=1}^5E^{5}_{p,1-p}$. By construction, we see that $E^{5}_{p,1-p}=Ker(\partial: G_1^{(1,\mathcal{A}_p))}\rightarrow G_0^{(1,\mathcal{A}_p)}))/G_1^{(1,\mathcal{A}_{p-1})}$. 
 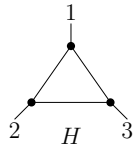
\begin{figure}[h]
\begin{center}
\scalebox{0.75}{
\begin{tikzpicture}
\draw[black](0,0)--(-0.7,-1);
\draw[black](0,0)--(0.7,-1);
\draw[black](0.7,-1)--(-0.7,-1);

\draw[black](0,0)--(0,0.4);
\draw[black](1,-1.3)--(0.7,-1);
\draw[black](-1,-1.3)--(-0.7,-1);

\fill[black] (0,0) circle (0.4ex);
\fill[black] (0.7,-1) circle (0.4ex);
\fill[black] (-0.7,-1) circle (0.4ex);
\node at (0,0.6){$1$};
\node at (-1,-1.5){$2$};
\node at (1,-1.5){$3$};
\node at (0,-1.6){$H$};
\end{tikzpicture}}
    \caption{The graph $H$ has no automorphisms, and there are six different possible orderings of its edges. This leads to six different generators $[H,\omega_{1,i}]$ in degree 1, for $i=1,...,6$, all related by the permutations of the ordering. }
    \label{fig:two-1eee}
    \end{center}
    \end{figure}

\noindent Now for every $\mathcal{A}_p$ the Kernel at the numerator has a single generator coming from the graph $H$ of Figure \ref{fig:two-1eee}, but since $H$ is $(1,\mathcal{A}_p)$-stable for any $p$ the generators coming from it belong to $G_1^{(1,\mathcal{A}_{p-1})}$ for any $p=2,...,5$, giving $E^{5}_{p,1-p}=0$. When $p=1$, $G_1^{(1,\mathcal{A}_{0})}=\{0\}$ by the convention we adopted at the beginning, so 
$$Gr^W_{6}H^{3}(\mathcal{M}_{1,3};\mathbb{Q})=\bigoplus_{p=1}^5E^{5}_{p,1-p}=E^{5}_{1,0}=Ker(\partial: G_1^{(1,\mathcal{A}_1))}\rightarrow G_0^{(1,\mathcal{A}_1)}))\cong \mathbb{Q}.$$ 
\end{itemize}
The computations made in this example give what was expected from Corollary 1.3 of \cite{CGP2}, which says that the top weight cohomology of $\mathcal{M}_{1,n}$ is supported in degree $n$ with rank $(n-1)!/2$ for $n\geq 3$, which equals 1 when $n=3$.
\end{myex}

\begin{myrem}
 We can use the Theorem to estimate the dimension of the cohomology of $\mathcal{
 M}_{g,n}$. Suppose we have $n,g\geq 1$, and a fixed sequence $[Ch_{\mathcal{A}_1}]\leq...\leq [Ch_{\mathcal{A}_p}]\leq...\leq [Ch_{\mathcal{A}_N}]$. Let $E^{\infty}_{p,k-p}$ be one of the pieces of the direct sum coming from the filtration, then $$dimH^{4g-6+2n-k}(\mathcal{M}_{g,n};\mathbb{Q})\geq dim Gr^W_{6g-6+2n}H^{4g-6+2n-k}(\mathcal{M}_{g,n};\mathbb{Q})\geq dimE^{\infty}_{p,k-p},$$ so if we are able to estimate the dimension of one of the pieces we can give nonvanishing results for the cohomology of $\mathcal{M}_{g,n}$.
\end{myrem}
\begin{myex}
Let $g=2$, $n=3$ and consider the sequence coming from the floor filtration of example \ref{floor}, extended on the left by the minimal chamber:
$$Ch_{\varepsilon^{(3)}}\leq Ch_{\mathcal{H}_3}\leq Ch_{1^{(3)}}.$$
By the Theorem \ref{LastTheorem} we know that  $Gr^W_{12}H^{8-k}(\mathcal{M}_{2,3};\mathbb{Q})\cong \bigoplus_{p=1}^NE^{\infty}_{p,k-p}.$ For $k=2$ the last term of the sum is $E^3_{3,-1}$, which has dimension at least $1$. To see this, consider the graphs in Figure \ref{finalexample}, then $H_1-H_2+H_3-G_1+G_2-G_3$ belongs to $G_2^{(2,3)}$ and one can see its differential is zero. However, its quotient by $G_2^{(2,\mathcal{H}_3)}$ is nonzero, namely it is $-G_1+G_2-G_3$, so it defines a nontrivial element of $E^3_{3,-5}$.

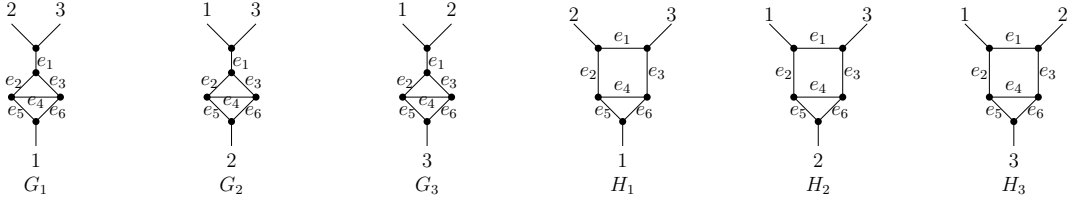
\begin{figure}[h]
\centering
\scalebox{0.65}{
\begin{tikzpicture}
\draw[black](0,0)--(-0.5,0.5);
\node at (0.2,-0.3){$e_1$};
\draw[black](0,0)--(0.5,0.5);
\fill[black] (0,0) circle (0.4ex);
\node at (-0.5,0.8){$2$};
\node at (0.5,0.8){$3$};
\draw[black](0,-0.5)--(0,0);
\fill[black] (0,-0.5) circle (0.4ex);
\draw[black](0,-0.5)--(0.5,-1);
\node at (0.45,-0.7){$e_3$};
\draw[black](0,-0.5)--(-0.5,-1);
\node at (-0.45,-0.7){$e_2$};
\fill[black] (-0.5,-1) circle (0.4ex);
\fill[black] (0.5,-1) circle (0.4ex);
\draw[black](-0.5,-1)--(0.5,-1);
\node at (0,-1.1){$e_4$};
\draw[black](-0.5,-1)--(0,-1.5);
\node at (-0.4,-1.3){$e_5$};
\draw[black](0.5,-1)--(0,-1.5);
\node at (0.45,-1.3){$e_6$};
\fill[black] (0,-1.5) circle (0.4ex);
\draw[black](0,-1.5)--(0,-2);
\node at (0,-2.3){$1$};
\node at (-0,-2.8){$G_1$};

\draw[black](0+4,0)--(-0.5+4,0.5);
\node at (0.2+4,-0.3){$e_1$};
\draw[black](0+4,0)--(0.5+4,0.5);
\fill[black] (0+4,0) circle (0.4ex);
\node at (-0.5+4,0.8){$1$};
\node at (0.5+4,0.8){$3$};
\draw[black](0+4,-0.5)--(0+4,0);
\fill[black] (0+4,-0.5) circle (0.4ex);
\draw[black](0+4,-0.5)--(0.5+4,-1);
\node at (0.45+4,-0.7){$e_3$};
\draw[black](0+4,-0.5)--(-0.5+4,-1);
\node at (-0.45+4,-0.7){$e_2$};
\fill[black] (-0.5+4,-1) circle (0.4ex);
\fill[black] (0.5+4,-1) circle (0.4ex);
\draw[black](-0.5+4,-1)--(0.5+4,-1);
\node at (0+4,-1.1){$e_4$};
\draw[black](-0.5+4,-1)--(0+4,-1.5);
\node at (-0.4+4,-1.3){$e_5$};
\draw[black](0.5+4,-1)--(0+4,-1.5);
\node at (0.45+4,-1.3){$e_6$};
\fill[black] (0+4,-1.5) circle (0.4ex);
\draw[black](0+4,-1.5)--(0+4,-2);
\node at (0+4,-2.3){$2$};
\node at (-0+4,-2.8){$G_2$};

\draw[black](0+8,0)--(-0.5+8,0.5);
\node at (0.2+8,-0.3){$e_1$};
\draw[black](0+8,0)--(0.5+8,0.5);
\fill[black] (0+8,0) circle (0.4ex);
\node at (-0.5+8,0.8){$1$};
\node at (0.5+8,0.8){$2$};
\draw[black](0+8,-0.5)--(0+8,0);
\fill[black] (0+8,-0.5) circle (0.4ex);
\draw[black](0+8,-0.5)--(0.5+8,-1);
\node at (0.45+8,-0.7){$e_3$};
\draw[black](0+8,-0.5)--(-0.5+8,-1);
\node at (-0.45+8,-0.7){$e_2$};
\fill[black] (-0.5+8,-1) circle (0.4ex);
\fill[black] (0.5+8,-1) circle (0.4ex);
\draw[black](-0.5+8,-1)--(0.5+8,-1);
\node at (0+8,-1.1){$e_4$};
\draw[black](-0.5+8,-1)--(0+8,-1.5);
\node at (-0.4+8,-1.3){$e_5$};
\draw[black](0.5+8,-1)--(0+8,-1.5);
\node at (0.45+8,-1.3){$e_6$};
\fill[black] (0+8,-1.5) circle (0.4ex);
\draw[black](0+8,-1.5)--(0+8,-2);
\node at (0+8,-2.3){$3$};
\node at (-0+8,-2.8){$G_3$};

\draw[black](0-0.5+12,0)--(1-0.5+12,0);
\draw[black](0-0.5+12,0)--(0-0.5+12,-1);
\draw[black](1-0.5+12,-1)--(1-0.5+12,0);
\draw[black](0-0.5+12,-1)--(1-0.5+12,-1);
\draw[black](0-0.5+12,-1)--(0.5-0.5+12,-1.5);
\draw[black](1-0.5+12,-1)--(0.5-0.5+12,-1.5);
\node at (0+12,0.2){$e_1$};
\node at (0+12,-0.8){$e_4$};
\node at (-0.7+12,-0.5){$e_2$};
\node at (0.7+12,-0.5){$e_3$};
\draw[black](0-0.5+12,0)--(-0.5-0.5+12,0.5);
\draw[black](1-0.5+12,0)--(1.5-0.5+12,0.5);
\draw[black](0.5-0.5+12,-2)--(0.5-0.5+12,-1.5);
\node at (-0.4+12,-1.3){$e_5$};
\node at (0.45+12,-1.3){$e_6$};
\node at (0+12,-2.3){$1$};
\node at (-1+12,0.7){$2$};
\node at (1+12,0.7){$3$};
\node at (-0+12,-2.8){$H_1$};
\fill[black] (0-0.5+12,0) circle (0.4ex);
\fill[black] (0.5+12,0) circle (0.4ex);
\fill[black] (0-0.5+12,-1) circle (0.4ex);
\fill[black] (0.5+12,-1) circle (0.4ex);
\fill[black] (0.5-0.5+12,-1.5) circle (0.4ex);

\draw[black](0-0.5+16,0)--(1-0.5+16,0);
\draw[black](0-0.5+16,0)--(0-0.5+16,-1);
\draw[black](1-0.5+16,-1)--(1-0.5+16,0);
\draw[black](0-0.5+16,-1)--(1-0.5+16,-1);
\draw[black](0-0.5+16,-1)--(0.5-0.5+16,-1.5);
\draw[black](1-0.5+16,-1)--(0.5-0.5+16,-1.5);
\node at (0+16,0.2){$e_1$};
\node at (0+16,-0.8){$e_4$};
\node at (-0.7+16,-0.5){$e_2$};
\node at (0.7+16,-0.5){$e_3$};
\draw[black](0-0.5+16,0)--(-0.5-0.5+16,0.5);
\draw[black](1-0.5+16,0)--(1.5-0.5+16,0.5);
\draw[black](0.5-0.5+16,-2)--(0.5-0.5+16,-1.5);
\node at (-0.4+16,-1.3){$e_5$};
\node at (0.45+16,-1.3){$e_6$};
\node at (0+16,-2.3){$2$};
\node at (-1+16,0.7){$1$};
\node at (1+16,0.7){$3$};
\node at (-0+16,-2.8){$H_2$};
\fill[black] (0-0.5+16,0) circle (0.4ex);
\fill[black] (0.5+16,0) circle (0.4ex);
\fill[black] (0-0.5+16,-1) circle (0.4ex);
\fill[black] (0.5+16,-1) circle (0.4ex);
\fill[black] (0.5-0.5+16,-1.5) circle (0.4ex);

\draw[black](0-0.5+20,0)--(1-0.5+20,0);
\draw[black](0-0.5+20,0)--(0-0.5+20,-1);
\draw[black](1-0.5+20,-1)--(1-0.5+20,0);
\draw[black](0-0.5+20,-1)--(1-0.5+20,-1);
\draw[black](0-0.5+20,-1)--(0.5-0.5+20,-1.5);
\draw[black](1-0.5+20,-1)--(0.5-0.5+20,-1.5);
\node at (0+20,0.2){$e_1$};
\node at (0+20,-0.8){$e_4$};
\node at (-0.7+20,-0.5){$e_2$};
\node at (0.7+20,-0.5){$e_3$};
\draw[black](0-0.5+20,0)--(-0.5-0.5+20,0.5);
\draw[black](1-0.5+20,0)--(1.5-0.5+20,0.5);
\draw[black](0.5-0.5+20,-2)--(0.5-0.5+20,-1.5);
\node at (-0.4+20,-1.3){$e_5$};
\node at (0.45+20,-1.3){$e_6$};
\node at (0+20,-2.3){$3$};
\node at (-1+20,0.7){$1$};
\node at (1+20,0.7){$2$};
\node at (-0+20,-2.8){$H_3$};
\fill[black] (0-0.5+20,0) circle (0.4ex);
\fill[black] (0.5+20,0) circle (0.4ex);
\fill[black] (0-0.5+20,-1) circle (0.4ex);
\fill[black] (0.5+20,-1) circle (0.4ex);
\fill[black] (0.5-0.5+20,-1.5) circle (0.4ex);
\end{tikzpicture}}
\caption{These are generators of $G_{2}^{(2,3)}$. The letters $e_i$ represent the chosen order on the set of edges. }
\label{finalexample}
\end{figure}
Then we conclude that $$dim H^{6}(\mathcal{M}_{2,3};\mathbb{Q})\geq dimGr^W_{12}H^{6}(\mathcal{M}_{2,3};\mathbb{Q})\geq dim(E^3_{3,-1})\geq 1,$$
i.e. $H^{6}(\mathcal{M}_{2,3};\mathbb{Q})$ does not vanish.

\end{myex}

\subsection*{Relative Homology} There is a relation between the spectral sequence associated to a filtration and the relative homology with respect to the inclusion.of $\Delta_{g,\mathcal{A}}\subset\Delta_{g,\mathcal{B}}$. As we already saw, whenever $X\subset Y$ is a subcomplex, for every $p\geq -1$ one can consider the exact sequence $$0 \rightarrow C_p(X)\rightarrow C_p(Y)\rightarrow C_p(Y,X)\rightarrow 0,$$
with $C_*(Y,X)$ computing the relative homology.

In particular, whenever we have an inclusion $\Delta_{g,\mathcal{A}}\hookrightarrow\Delta_{g,\mathcal{B}}$, the relative chain complex $C_*(\Delta_{g,\mathcal{B}},\Delta_{g,\mathcal{A}})$ has its homology coinciding with relative rational homology
$$H_i(C_*(\Delta_{g,\mathcal{B}},\Delta_{g,\mathcal{A}}))\cong H_i(\Delta_{g,\mathcal{B}},\Delta_{g,\mathcal{A}};\mathbb{Q}),$$ as it has a natural isomorphism with $C_*(\Delta_{g,\mathcal{B}})/C_*(\Delta_{g,\mathcal{A}})$.

Now as we saw in the proof of Theorem \ref{HomTh} we have a decomposition $C_*(\Delta_{g,\mathcal{A}})=A^{(g,\mathcal{A})}\bigoplus B^{(g,\mathcal{A})}$, so we can decompose also $C_*(\Delta_{g,\mathcal{B}})/C_*(\Delta_{g,\mathcal{A}})$ into $A^{(g,\mathcal{B})}/C_*(\Delta_{g,\mathcal{A}})\bigoplus B^{(g,\mathcal{B})}/C_*(\Delta_{g,\mathcal{A}})$. 

The term $B^{(g,\mathcal{B})}/C_*(\Delta_{g,\mathcal{A}})$ is acyclic, since also $B^{(g,\mathcal{B})}$ is. Now consider the shifting degree injection $j:G^{(g,\mathcal{A})}\rightarrow C_*(\Delta_{g,\mathcal{A}},\mathbb{Q})$ of Theorem \ref{HomTh}: it is clear by the definition of $A^{(g,\mathcal{B})}$ that $A^{(g,\mathcal{B})}/C_*(\Delta_{g,\mathcal{A}})$ is isomorphic to $A^{(g,\mathcal{B})}/j(G^{(g,\mathcal{A})})$ and this is isomorphic to $G^{(g,\mathcal{B})}/G^{(g,\mathcal{A})}:=G^{(g,\mathcal{B},\mathcal{A})}$, with the isomorphism shifting degrees by $2g-1$. The complex $G^{(g,\mathcal{B},\mathcal{A})}$ can be seen as the one generated by $(g,\mathcal{B})$-stable but not $(g,\mathcal{A})$-stable graphs, with the same conventions on the degree and the same boundary map.

Through the latter isomoprhism, we can conclude that there is an isomorphism $$H_{k-2g+1}(G^{(g,\mathcal{B},\mathcal{A})})\cong H_k(\Delta_{g,\mathcal{B}},\Delta_{g,\mathcal{A}};\mathbb{Q}).$$

\begin{myex}
Consider the floor filtration of Example \ref{floor}. A graph is $(g,\mathcal{H}_l)$-stable if and only if its leaves have at least $l$ markings, and its vertices of valence 2 have at least one. In particular, when $l=2$ we get the usual notion of stability. Then the complex $G^{(g,\mathcal{H}_l,\mathcal{H}_{l+i})}$ is generated by stable graphs (in the standard sense) with a number of markings on each leaf between $l$ and $l-i-1$., and by the previous computations we have $H_{k-2g+1}(G^{(g,\mathcal{H}_l,\mathcal{H}_{l+i})})\cong H_k(\Delta_{g,\mathcal{H}_l},\Delta_{g,\mathcal{H}_{l+i}};\mathbb{Q}).$ 

In particular, if $i=1$, $G^{(g,\mathcal{H}_l,\mathcal{H}_{l+1})}$ is generated by stable graphs with exactly $l$ markings on each leaf. 

When $l=2$, we have $\mathcal{H}_2\in Ch_{1^{(n)}}$ so $$H_{k-2g+1}(G^{(g,\mathcal{H}_2,\mathcal{H}_{3})})\cong H_k(\Delta_{g,n},\Delta_{g,\mathcal{H}_{3}};\mathbb{Q}).$$

As an example of computation, if $g=1$ and for every $n$ we have $H_0(\Delta_{1,n},\Delta_{1,\mathcal{H}_3};\mathbb{Q})=0$, $H_1(\Delta_{1,n},\Delta_{1,\mathcal{H}_3};\mathbb{Q})=0$,
\end{myex}

Now for a given filtration, at each step we have $E^1_{p,q}=H_{p+q}(G^{(g,\mathcal{A}_p,\mathcal{A}_{p-1})})$, so we conclude that the first page of the spectral sequence computes the relative homology at each step of the filtration.

\nocite{*}
\bibliographystyle{alpha}
\bibliography{Bibliography}

\end{document}